\documentclass[a4paper]{amsart}
\usepackage{amsthm,amsmath,amssymb,amscd,multirow}
\usepackage[latin1]{inputenc}
 \newtheorem{theorem}{Theorem}
 \newtheorem{main}{Theorem}
 \newtheorem{proposition}[theorem]{Proposition}
 \newtheorem{lemma}[theorem]{Lemma}

 \theoremstyle{definition}
 \newtheorem{definition}[theorem]{Definition}

 \theoremstyle{remark}
 \newtheorem{remark}[theorem]{Remark}
 
 \numberwithin{equation}{section}
 \numberwithin{theorem}{section}

 \newcommand{\Spec}{{\rm Spec}}
 \newcommand{\Spf}{{\rm Spf}}
 \newcommand{\vol}{{\rm vol}}
 \newcommand{\Quot}{{\rm Quot}}
 \newcommand{\Stab}{{\rm Stab}}
 \newcommand{\End}{{\rm End}}
 \newcommand{\Nilp}{{\rm Nilp}}
 \newcommand{\Sets}{{\rm Sets}}
 \newcommand{\hei}{{\rm ht}}
 \newcommand{\id}{{\rm id}}

 \newcommand{\et}{{\rm et}}
 \newcommand{\m}{{\rm m}}

\newcommand{\np}{{\rm np}}
 \newcommand{\bi}{{\rm bi}}
 \newcommand{\ad}{{\rm ad}}
 \newcommand{\rk}{{\rm rk}}
 
 \newcommand{\Ann}{{\rm Ann}}
 \newcommand{\Hom}{{\rm Hom}}
 \newcommand{\red}{{\rm red}}
 
 \newcommand{\defe}{{\rm def}}
 
 \newcommand{\F}{{\boldsymbol F}}

\newcommand{\LT}{{\rm LT}}

\begin{document}
\begin{title}
{The global structure of moduli spaces of polarized $p$-divisible groups}
\end{title}
\author{Eva Viehmann}
\address{Mathematisches Institut der Universit\"{a}t Bonn\\ Beringstrasse 1\\53115 Bonn\\Germany}
\date{}
\begin{abstract}
We study the global structure of moduli spaces of
quasi-isogenies of polarized $p$-divisible groups introduced by Rapoport and Zink. Using the corresponding results for non-polarized $p$-divisible groups from a previous paper, we determine their dimensions and their sets of connected components and of irreducible components.
\end{abstract}
\maketitle 

\section{Introduction}

Let $k$ be an algebraically closed field of characteristic $p>2$. Let $W=W(k)$ be its ring of Witt vectors and $L=\Quot (W)$. Let $\sigma$ be the Frobenius automorphism on $k$ as well as on $W$. By $\Nilp_W$ we denote the category of schemes $S$ over $\Spec(W)$ such that $p$ is locally nilpotent on $S$. Let $\overline{S}$ be the closed subscheme of $S$ that is defined by the ideal sheaf $p\mathcal{O}_S$. Let $(\mathbb{X},\lambda_{\mathbb{X}})$ be a principally polarized $p$-divisible group over $k$. If $X$ is a $p$-divisible group, we denote its dual by $\hat{X}$. Then the polarization $\lambda_{\mathbb{X}}$ is an isomorphism $\mathbb{X}\rightarrow \hat{\mathbb{X}}$.

We consider the functor
$$\mathcal{M}:\Nilp_W\rightarrow  \Sets,$$ which
assigns to $S\in \Nilp_W$ the set of isomorphism classes of pairs
$(X,\rho)$, where $X$ is a $p$-divisible group over $S$ and $\rho:
\mathbb{X}_{\overline{S}}=\mathbb{X}\times_{\Spec
(k)}\overline{S}\rightarrow X\times_{S}\overline{S}$ is a
quasi-isogeny such that the following condition holds. There exists a quasi-isogeny $\lambda:X\rightarrow \hat X$ such that $\rho^{\vee}\circ\lambda_{\overline{S}}\circ\rho$ and $\lambda_{\mathbb{X},\overline{S}}$ coincide up to a scalar. Two pairs $(X_1, \rho_1)$ and $(X_2, \rho_2)$ are
isomorphic if $\rho_1\circ\rho_2^{-1}$ lifts to an isomorphism
$X_2\rightarrow X_1$. This functor is representable by a formal
scheme $\mathcal{M}$ which is locally formally of finite type
over $\Spf(W)$ (see \cite{RapoportZink}, Thm. 3.25). Let
$\mathcal{M}_{\red}$ be its underlying reduced subscheme, that is the reduced subscheme of $\mathcal{M}$ defined by the maximal ideal of definition. Then $\mathcal{M}_{\red}$ is a scheme over $\Spec (k)$.

The analogues of these moduli spaces for $p$-divisible groups without polarization have been studied by de Jong and Oort in \cite{deJongOort} for the case that the rational Dieudonn\'{e} module of $\mathbb{X}$ is simple and in \cite{modpdiv} without making this additional assumption. There, the sets of connected components and of irreducible components, as well as the dimensions, are determined. In the polarized case, the moduli spaces $\mathcal{M}_{\red}$ have been examined in several low-dimensional cases. For example, Kaiser (\cite{Kaiser}) proves a twisted fundamental lemma for $GSp_4$ by giving a complete description in the case that $\mathbb{X}$ is two-dimensional and supersingular. An independent description of this case is given by Kudla and Rapoport in \cite{KudlaRapoport}. In \cite{Richartz}, Richartz describes the moduli space in the case of three-dimensional supersingular $\mathbb{X}$. In this paper we derive corresponding results on the global structure of the moduli space $\mathcal{M}_{\red}$ for arbitrary $\mathbb{X}$. 

The first main result of this paper concerns the set of connected components of $\mathcal{M}_{\red}$.

\begin{main}\label{thmconncomp} Let $\mathbb{X}$ be nontrivial and let $\mathbb{
X}_{\m}\times \mathbb{X}_{\bi}\times \mathbb{X}_{\et}$ be the decomposition into its multiplicative, bi-infinitesimal, and \'{e}tale part. Then $$\pi_0(\mathcal{M}_{\red})\cong
 \left(GL_{\hei(\mathbb{X}_{\m})}(\mathbb{Q}_p)/GL_{\hei(\mathbb{X}_{\m})}(\mathbb{Z}_p)\right)\times \mathbb{Z}.$$
\end{main}

Next we consider the set of irreducible components of
$\mathcal{M}_{\red}$. Let $(N,F)$ be the rational Dieudonn\'{e} module of $\mathbb{X}$. Here, $N$ is an $L$-vector space of dimension $\hei(\mathbb{X})$ and $F:N\rightarrow N$ is a $\sigma$-linear isomorphism. The polarization $\lambda_{\mathbb{X}}$ induces an anti-symmetric bilinear perfect pairing $\langle\cdot,\cdot\rangle$ on $N$. Let $G$ be the corresponding general symplectic group of automorphisms of $N$ respecting $\langle\cdot,\cdot\rangle$ up to a scalar. Let $$J=\{g\in G(L) \mid g\circ F=F\circ g\}.$$ It is the set of $\mathbb{Q}_p$-valued points of an algebraic group over $\mathbb{Q}_p$ (see \cite{RapoportZink}, Prop. 1.12). There is an action of $J$ on $\mathcal{M}_{\red}$. 

\begin{main}\label{thmirrkomp}
The action of $J$ on the set of irreducible
components of $\mathcal{M}_{\red}$ is transitive.
\end{main}

We choose a decomposition $N=\bigoplus_{j=1}^{l}N^{j}$ with $N^{j}$ simple of slope $\lambda_j=m_j/(m_j+n_j)$ with $(m_j,n_j)=1$ and
$\lambda_j\leq\lambda_{j'}$ for $j<j'$. Let 
\begin{equation*}m=\left\lfloor\frac{1}{2}\sum_j \min\{m_j,n_j\}\right\rfloor,\end{equation*} where $\lfloor x \rfloor$ is the greatest integer less or equal $x$. As $N$ is the isocrystal of a polarized $p$-divisible group, its Newton polygon is symmetric, i.~e.~ $\lambda_{l+1-j}=1-\lambda_j$. Hence we obtain
\begin{equation}\label{gldefm}m=\left\lfloor\sum_{\{j\mid m_j<n_j\}} m_j+\frac{1}{2}\mid\{j\mid m_j=n_j=1\}\mid\right\rfloor.\end{equation}

\begin{main}\label{thmdim}
$\mathcal{M}_{\red}$ is equidimensional of dimension
\begin{equation}\label{gldimform}
\dim\mathcal{M}_{\red}=\frac{1}{2}\left(\sum_{j}\frac{(m_j-1)(n_j-1)}{2}+\sum_{j<j'}m_jn_{j'}+ m\right).
\end{equation}
\end{main}

Note that the equidimensionality  is already a consequence of Theorem \ref{thmirrkomp}. However, it also follows from the proof of the dimension formula without requiring additional work.

Our results on the set of connected components and on the dimension of $\mathcal{M}_{\red}$ are analogous to those for other affine Deligne-Lusztig sets for split groups where a scheme structure is known. We now define these affine Deligne-Lusztig varieties and give a brief overview over the general results in comparison to the results for the case treated in this paper.

Let  $\mathcal{O}$ be a finite extension of $\mathbb{Z}_p$ or $\mathbb{F}_p[[t]]$ and let $G$ be a split connected reductive group over $\mathcal{O}$. Let $\F$ be the quotient field of $\mathcal{O}$. Let $K=G(\mathcal{O})$. Let $L$ be the completion of the maximal unramified extension of $\F$ and let $\sigma$ be the Frobenius of $L$ over $\F$. Let $A$ be a maximal torus and $B$ a Borel subgroup containing $A$. Let $\mu\in X_*(A)$ be dominant and let $b\in G(L)$ be in $B(G,\mu)$. This last condition ensures that the affine Deligne-Lusztig variety will be nonempty (compare \cite{Rapoport}, 5). Let $\varepsilon^{\mu}$ be the image of $p$ or $t\in \F^{\times}$ under $\mu$. Let 
\begin{equation}
X_{\mu}(b)=\{g\in G(L)/K\mid g^{-1}b\sigma(g)\in K\varepsilon^{\mu} K\}
\end{equation}
be the generalized affine Deligne-Lusztig set associated to $\mu$ and $b$. There are two cases where it is known that $X_{\mu}(b)$ is the set of $k$-valued points of a scheme. Here, $k$ denotes the residue field of $\mathcal{O}_L$. The first case is that $\F=\mathbb{Q}_p$ and that $X_{\mu}(b)$ is the set of $k$-valued points of a Rapoport-Zink space of type (EL) or (PEL). Let for example $G=GL_h$ or $GSp_{2h}$. The case $G=GL_h$ was considered in \cite{modpdiv}. For the moduli spaces considered in this paper choose a basis  $\{e_i,f_i\mid 1\leq i\leq h\}$ identifying $N$ with $L^{2h}$ and the symplectic form on $N$ with the symplectic form on $L^{2h}$ defined by requiring that $\langle e_i,e_j\rangle=\langle f_i,f_j\rangle=0$ and $\langle e_i,f_j\rangle=\delta_{i,h+1-j}$. Then let $G=GSp_{2h}$. Let $B$ be the Borel subgroup of $G$ fixing the complete isotropic flag $(e_1)\subset(e_1,e_2)\subset \dotsm\subset (e_1,\dotsc, e_h)$. We choose $A$ to be the diagonal torus. Let $\pi_1(G)$ be the quotient of $X_*(A)$ by the coroot lattice of $G$. Then $\pi_1(G)\cong \mathbb{Z}$, where $\alpha\in X_*(A)$ is mapped to $v_p(c)$ with $c$ being the multiplier of $p^{\alpha}\in GSp_{2h}(L)$. Let $\mu\in X_*(A)$ be the unique minuscule element whose image in $\pi_1(G)$ is 1. Then $p^{\mu}$ is a diagonal matrix with diagonal entries $1$ and $p$, each with multiplicity $h$. We write $F=b\sigma$ with $b\in G$.  Note that there is a bijection between $\mathcal{M}_{\red}(k)$ and the set of Dieudonn\'{e} lattices in $N$. Using the above notation, we have the bijection
\begin{eqnarray*}
X_{\mu}(b)&\rightarrow& \mathcal{M}_{\red}(k)\\
g&\mapsto&g(W(k)^{2h}).
\end{eqnarray*}
The second case is that $\F$ is a function field. Here $X_{\mu}(b)$ obtains its scheme structure by considering it as a subset of the affine Grassmannian $G(L)/G(\mathcal{O}_L)$. In this case we do not have to assume $\mu$ to be minuscule. The $X_{\mu}(b)$ are locally closed subschemes of the affine Grassmannian. The closed affine Deligne-Lusztig varieties $X_{\preceq\mu}(b)$ are defined to be the closed reduced subschemes of $G(L)/G(\mathcal{O}_L)$ given by $X_{\preceq\mu}(b)=\bigcup_{\mu'\preceq\mu} X_{\mu'}(b)$. Here $\mu'\preceq\mu$ if $\mu-\mu'$ is a non-negative linear combination of positive coroots. Note that the two schemes $X_{\mu}(b)$ and $X_{\preceq\mu}(b)$ coincide if $\mu$ is minuscule.

The sets of connected components of the moduli spaces of non-polarized $p$-divisible groups are given by a formula completely analogous to Theorem \ref{thmconncomp} (compare \cite{modpdiv}, Thm. A). For closed affine Deligne-Lusztig varieties in the function field case, the set of connected components is also given by a generalization of the formula in Theorem \ref{thmconncomp} (see \cite{conncomp}, Thm. 1). The sets of connected components of the non-closed $X_{\mu}(b)$ are not known in general. There are examples (compare \cite{conncomp}, Section 3) which show that a result analogous to Theorem \ref{thmconncomp} cannot hold for all non-closed $X_{\mu}(b)$.

The only further general case where the set of irreducible components is known are the reduced subspaces of moduli spaces of $p$-divisible groups without polarization. Here, the group $J$ also acts transitively on the set of irreducible components. There are examples of affine Deligne-Lusztig varieties in the function field case associated to non-minuscule $\mu$ where this is no longer true (compare \cite{dimdlv}, Ex. 6.2).

To discuss the formula for the dimension let us first reformulate Theorem \ref{thmdim}. Let $G=GSp_{2h}$ and $\mu$ be as above. Let $\nu=(\lambda_i)\in \mathbb{Q}^h\cong X_*(A)_{\mathbb{Q}}$ be the (dominant) Newton vector associated to $(N,F)$ as defined by Kottwitz, see \cite{Kottwitz1}. Let $\rho$ be the half-sum of the positive roots of $G$ and $\omega_i$ the fundamental weights of the adjoint group $G_{\ad}$. Then one can reformulate (\ref{gldimform}) as
\begin{equation}\label{gldimform3}
\dim \mathcal{M}_{\red}=\langle 2\rho,\mu-\nu\rangle +\sum_i\lfloor\langle \omega_i,\nu-\mu\rangle\rfloor.
\end{equation}
In this form, the dimension formula proves a special case of a conjecture by Rapoport (see \cite{Rapoport}, Conjecture 5.10) for the dimension of affine Deligne-Lusztig varieties. Denote by $\rk_{\mathbb{Q}_p}$ the dimension of a maximal $\mathbb{Q}_p$-split subtorus and let $\defe_G(F)=\rk_{\mathbb{Q}_p}G-\rk_{\mathbb{Q}_p}J$. In our case, this is equal to $h-\lceil l/2\rceil$ where $l$ is the number of simple summands of $N$. Using Kottwitz's reformulation of the right hand side of (\ref{gldimform3}) in \cite{Kottwitz2}, we obtain
\begin{equation}\label{gldimform2}
\dim \mathcal{M}_{\red}=\langle \rho,\mu-\nu\rangle -\frac{1}{2}\defe_G(F).
\end{equation}

For the case of moduli of $p$-divisible groups for $G=GL_h$, the analogous formula for the dimension is shown in \cite{modpdiv}. In the function field case, the dimension of the generalized affine Deligne-Lusztig variety has been determined in \cite{dimdlv}, \cite{GHKR}. The formula for the dimension is also in this case the analogue of (\ref{gldimform2}). 

The dimension of the moduli spaces $\mathcal{M}_{\red}$ is also studied by Chai and Oort using a different approach. In \cite{Oort3}, Oort defines an almost product structure (that is, up to a finite morphism) on Newton strata of moduli spaces of polarized abelian varieties. It is given by an isogeny leaf and a central leaf for the $p$-divisible group. The dimension of the isogeny leaf is the same as that of the corresponding $\mathcal{M}_{\red}$. He announces a joint paper with Chai, in which they prove a dimension formula for central leaves (compare \cite{Oort3}, Remark 2.8). The dimension of the Newton polygon stratum itself is known from \cite{Oort1}. Then the dimension of $\mathcal{M}_{\red}$ can also be computed as the difference of the dimensions of the Newton polygon stratum and the central leaf. 

We outline the content of the different sections of the paper. In Section \ref{sec2} we introduce the necessary background and notation, and reduce the problem to the case of bi-infinitesimal groups. In the third and fourth section, we define the open dense subscheme $\mathcal{S}_1$ where the $a$-invariant of the $p$-divisible group is $1$ and describe its set of closed points. This description is refined in Sections 5 and 6 to prove the theorems on the set of irreducible components and on the dimension, respectively. In the last section we determine the set of connected components.

\noindent{\it Acknowledgement.} Part of this paper was written during a stay at the Universit\'{e} Paris-Sud at Orsay which was supported by a fellowship within the Postdoc-Program of the German Academic Exchange Service (DAAD). I thank the Universit\'{e} Paris-Sud for its hospitality. 

\section{Notation and preliminary reductions}\label{sec2}

\subsection{A decomposition of the rational Dieudonn\'{e} module}
The principal polarization $\lambda_{\mathbb{X}}$ equips the rational Dieudonn\'{e} module $(N,F)$ of $\mathbb{X}$ with a nondegenerate anti-symmetric bilinear pairing $\langle \cdot,\cdot\rangle$. It satisfies
\begin{equation}\label{glskpr}
\langle v,Fw\rangle=\sigma(\langle Vv,w\rangle)
\end{equation}
for all $x,y\in N$.

We assumed $k$ to be algebraically closed. Then an easy calculation using the classification of isocrystals shows that $N$ has a decomposition into subisocrystals $N_i$ of one of the following types. Let $l$ be the number of supersingular summands in a decomposition of $N$ into simple isocrystals. Then
\begin{equation}
N=\begin{cases}
N_0\oplus N_1& \text{if }l\text{ is even}\\
N_0\oplus N_{\frac{1}{2}}\oplus N_1&\text{otherwise.}
\end{cases}
\end{equation}
Here the slopes of $N_0$ are smaller or equal to $\frac{1}{2}$ and $N_1$ is the isocrystal dual to $N_0$. The summand $N_{\frac{1}{2}}$ is simple and supersingular. Especially,
$$\langle N_0,N_0\rangle=\langle N_1,N_1\rangle=\langle N_0,N_{\frac{1}{2}}\rangle=\langle N_1,N_{\frac{1}{2}}\rangle=0.$$
For $i\in\{0,\frac{1}{2},1\}$ we denote by $p_i$ the canonical projection $N\rightarrow N_i$.

The moduli spaces $\mathcal{M}_{\red}$ for different $(\mathbb{X},\lambda_{\mathbb{X}})$ in the same isogeny class are isomorphic. Replacing $\mathbb{X}$ by an isogenous group we may assume that 
\begin{equation}
\mathbb{X}=\begin{cases}
\mathbb{X}_0\times\mathbb{X}_1& \text{if }l\text{ is even}\\
\mathbb{X}_0\times \mathbb{X}_{\frac{1}{2}}\times \mathbb{X}_1&\text{otherwise.}
\end{cases}
\end{equation}
Here, $\mathbb{X}_i$ is such that its rational Dieudonn\'{e} module is $N_i$.

Mapping $(X,\rho)\in \mathcal{M}_{\red}(k)$ to the Dieudonn\'{e} module of $X$ defines a bijection between $\mathcal{M}_{\red}(k)$ and the set of Dieudonn\'{e} lattices in $N$ that are self-dual up to a scalar. Here a sublattice $\Lambda$ of $N$ is called a Dieudonn\'{e} lattice if $\varphi(\Lambda)\subseteq \Lambda$ for all $\varphi$ in the Dieudonn\'{e} ring of $k$,
\begin{equation}\label{gldefd}
\mathcal{D}=\mathcal{D}(k)=W(k)[F,V]/(FV=VF=p,aV=V\sigma(a),Fa=\sigma(a)F).
\end{equation} 
All lattices considered in this paper are Dieudonn\'{e} lattices. A lattice $\Lambda\subset N$ is self-dual up to a scalar if the dual lattice $\Lambda^{\vee}$ satisfies $\Lambda^{\vee}=c\Lambda$ with $c\in L^{\times}$.

The following notion is introduced by Oort in \cite{Oort}.
\begin{definition}\label{defminimal}
Let $X$ be a $p$-divisible group over $k$ and $\Lambda_{\min}$ its Dieudonn\'{e} module. Assume that $\Lambda_{\min}$ is the direct sum of submodules $\Lambda_{\min}^i$ such that $N^i=\Lambda_{\min}^i\otimes_W L$ is simple and that $\End (\Lambda_{\min}^i)$ is a maximal order in $\End (N^i)$. Then $X$ is called a minimal $p$-divisible group.
\end{definition}
\begin{lemma}\label{lemminimal}
There is a $k$-valued point $(X,\rho)$ of $\mathcal{M}_{\red}$ such that $X$ is minimal.
\end{lemma}
\begin{proof}
Let $N_0$ and $N_1$ as in the decomposition above. Let $\Lambda_{\min,0}\subset N_0$ be the lattice of a minimal $p$-divisible group and let $\Lambda_{\min,\frac{1}{2}}\subset N_{\frac{1}{2}}$ be the Dieudonn\'{e} module of $\mathbb{X}_{\frac{1}{2}}$. There is only one isomorphism class of one-dimensional supersingular $p$-divisible groups and it consists of minimal $p$-divisible groups. Let $c\in L^{\times}$ with $\Lambda_{\min,\frac{1}{2}}^{\vee}=c \Lambda_{\min,\frac{1}{2}}$. Let $$\Lambda_{\min,1}=\{x\in N_1\mid \langle x,cy\rangle\in W\text{ for all }y\in \Lambda_{\min,0}\}.$$ Then $\Lambda_{\min,1}$ is also the Dieudonn\'{e} module of a minimal $p$-divisible group. Furthermore, $\Lambda_{\min}=\Lambda_{\min,0}\oplus \Lambda_{\min,\frac{1}{2}}\oplus \Lambda_{\min,1}$ satisfies $\Lambda_{\min}^{\vee}=c \Lambda_{\min}$. Thus $\Lambda_{\min}$ corresponds to an element of $\mathcal{M}_{\red}(k)$ and to a minimal $p$-divisible group.
\end{proof}
\begin{remark}\label{remmin}
There is the following explicit description of the Dieudonn\'{e} module of a minimal $p$-divisible group: For each $j$, there is a basis $e_1^j,\dotsc,e_{m_j+n_j}^j$ of $N^j$ with $F(e_i^j)=e_{i+m_j}^j$ for all $i,j$. Here we use the notation $e^j_{i+m_j+n_j}=pe_i^j$. For the existence compare for example \cite{modpdiv}, 4.1. Furthermore, these bases may be chosen such that $\langle e_i^j,e_{i'}^{j'}\rangle=\delta_{j,l+1-j'}\cdot\delta_{i,m_j+n_j+1-i'}$ for $1\leq i,i'\leq m_j+n_j=m_{l+1-j}+n_{l+1-j}$. Then we can take the lattice $\Lambda_{\min}$ to be the lattice generated by these basis elements $e_i^j$. Using these bases to identify $N$ with $L^{2h}$ has the advantage that $\Stab(\Lambda_{\min})=K$.
\end{remark}

\subsection{Moduli of non-polarized $p$-divisible groups}\label{secnp}

For the moment let $\mathbb{X}$ be a $p$-divisible group without polarization. Then associated to $\mathbb{X}$ there is an analogous moduli problem of quasi-isogenies of $p$-divisible groups without polarization. If $\mathbb{X}$ is polarized, we thus obtain two functors which are closely related. In this section we recall the definition of the moduli spaces of non-polarized $p$-divisible groups and relate them to $\mathcal{M}_{\red}$. Besides, we provide a technical result on lattices in isocrystals which we need in the following section. 

Let $\mathcal{M}_{\mathbb{X}}^{\np}$ be the functor associating to a scheme $S\in \Nilp_W$ the set of pairs $(X,\rho)$ where $X$ is a $p$-divisible group over $S$ and $\rho$ a quasi-isogeny $\mathbb{X}_{\overline{S}}\rightarrow X_{\overline{S}}$. Two such pairs $(X_1,\rho_1)$ and $(X_2,\rho_2)$ are identified in this set if $\rho_1\circ\rho_2^{-1}$ lifts to an isomorphism $X_2\rightarrow X_1$ over $S$. This functor is representable by a formal scheme which is locally formally of finite type over $\Spf(W)$ (see \cite{RapoportZink}, Theorem 2.16). Let $\mathcal{M}_{\mathbb{X},\red}^{\np}$ be its reduced subscheme. We always include $\mathbb{X}$ in this notation, because we compare $\mathcal{M}_{\red}$ to the two moduli spaces $\mathcal{M}_{\mathbb{X},\red}^{\np}$ and $\mathcal{M}_{\mathbb{X}_0,\red}^{\np}$.

Let $J^{\np}=\{g\in GL(N)\mid g\circ F=F\circ g\}$. Then $J\subseteq J^{\np}$.

If $\mathbb{X}$ is a principally polarized $p$-divisible group, then forgetting the polarization induces a natural inclusion as a closed subscheme
\begin{equation*}
\mathcal{M}_{\red}\hookrightarrow\mathcal{M}_{\mathbb{X},\red}^{\np}.
\end{equation*}

Furthermore, there is a natural inclusion as a closed subscheme
\begin{equation}\label{glincx0}
\mathcal{M}_{\mathbb{X}_0,\red}^{\np}\hookrightarrow\mathcal{M}_{\red}
\end{equation}
mapping $(X_0,\rho_0)$ to $(X_0\times X_0^{\vee},(\rho,\rho^{\vee}))$ if the number of supersingular summands of $N$ is even and to $(X_0\times X_{\frac{1}{2}}\times X_0^{\vee},(\rho,\rho_{\frac{1}{2}},\rho^{\vee}))$ otherwise. Here $X_{\frac{1}{2}}=\mathbb{X}_{\frac{1}{2}}$ is the unique one-dimensional supersingular  $p$-divisible group and $\rho_{\frac{1}{2}}=\id$.

Let $\tilde{v}$ be the valuation on the Dieudonn\'{e} ring $\mathcal{D}$ determined by 
\begin{equation}\label{gldefvtilde}
\tilde{v}(aF^iV^j)=2v_p(a)+i+j
\end{equation}
for every $a\in W(k)$. One can decompose each $B\in\mathcal{D}$ uniquely as $B=\LT(B)+B'$ with $\tilde{v}(B')>\tilde{v}(B)$ and $$\LT(B)=\sum_{0\leq i\leq \tilde{v}(B),2\alpha+i=\tilde{v}(B)}p^{\alpha}([a_i]V^i+[b_i]F^i).$$ Here $[a_i]$ and $[b_i]$ are Teichm\"{u}ller representatives of elements of $k^{\times}$ or 0 and $[b_0]=0$.

\begin{lemma}\label{lemdefa}
Let $(N_0,b_0\sigma)$ be the rational Dieudonn\'{e} module of some $p$-divisible group over $k$. Let $m=v_p(\det b_0)$ and $n=\dim_L (N_0)-m$. Let $v\in N_0$ be not contained in any proper sub-isocrystal of $N_0$. 
\begin{enumerate}
\item $\Ann(v)=\{\varphi\in\mathcal{D}\mid \varphi(v)=0\}$ is a principal left ideal of $\mathcal{D}$. There is a generating element of the form $$A=aF^n+bV^m+\sum_{i=0}^{n-1}a_iF^i+\sum_{i=1}^{m-1}b_iV^i.$$
with $a,b\in W^{\times}$ and $a_i,b_i\in W$. 
\item 
If $N_0$ is simple (and thus of slope $m/(m+n)$), we have 
\begin{equation}\label{gllt}
\LT(A)=\begin{cases}
[a]F^n&\text{if }n<m\\
[b]V^m&\text{if }m<n\\
[a]F+[b]V&\text{if }m=n=1
\end{cases}
\end{equation}
for some $a,b\in  k^{\times} $.
\item Let $N_0=\oplus_j N^j$ be a decomposition of $N_0$ into simple summands. Then $\LT(A)=\LT(\prod_{j}L_j).$ Here each $L_j$ is of the form (\ref{gllt}) associated to some nonzero element in $N^j$.
\end{enumerate}
\end{lemma}
\begin{proof}
We use induction on the number of summands in a decomposition of $N_0$ into simple isocrystals. If $N_0$ is simple, the lemma follows immediately from \cite{modpdiv}, Lemma 4.12. For the induction step write $N_0=N'\oplus N^{''}$ where $N'$ is simple. Let $A'$ be as in the lemma and associated to $N'$ and $p_{N'}(v)$. Note that an element of an isocrystal is not contained in any proper sub-isocrystal if and only if the Dieudonn\'{e} module generated by the element is a lattice. Let $\Lambda$ be the lattice generated by $v$. The Dieudonn\'{e} module generated by $A'(v)$ is equal to $\Lambda\cap N^{''}$, and hence also a lattice. We may therefore apply the induction hypothesis to $A'(v)$ and $N^{''}$ and obtain some $A^{''}$ generating $\Ann(A'(v))$. Thus $\Ann(v)$ is a principal left ideal generated by $A^{''}A'$. Multiplying the corresponding expressions for $A^{''}$ and $A'$, the lemma follows.
\end{proof}

\subsection{Reduction to the bi-infinitesimal case}\label{secredbi}
Let $\mathbb{X}=\mathbb{X}_{\et}\times\mathbb{X}_{\bi}\times\mathbb{X}_{\m}$ be the decomposition of $\mathbb{X}$ into its \'{e}tale, bi-infinitesimal, and multiplicative parts.
\begin{lemma} We have
$$\mathcal{M}_{\mathbb{X},\red}\cong \begin{cases}\mathcal{M}^{\np}_{\mathbb{X}_{\et},\red}\times\mathcal{M}_{\mathbb{X}_{\bi},\red}&\text{if }\mathbb{X}_{\bi}\text{ is nontrivial}\\
\mathcal{M}^{\np}_{\mathbb{X}_{\et},\red}\times\mathbb{Z}&\text{otherwise}.
\end{cases}$$ and $$\mathcal{M}^{\np}_{\mathbb{X}_{\et},\red}\cong GL_{\hei(\mathbb{X}_{\et})}(\mathbb{Q}_p)/GL_{\hei(\mathbb{X}_{\et})}(\mathbb{Z}_p).$$
\end{lemma}
\begin{proof}
Consider the following morphism from the right to the left hand side of the first isomorphism. In the first case, a point $((X_{\et},\rho_{\et}),( X_{\bi},\rho_{\bi}))$ is mapped to $(X_{\et}\times X_{\bi}\times X_{\m},(\rho_{\et},\rho_{\bi},\rho_{\m}))$ where $X_{\m}=\hat{X}_{\et}$. Furthermore, $\rho_{\m}$ is the dual isogeny of $p^l\cdot\rho_{\et}$ and $p^l$ is the scalar determined by the duality condition for $\rho_{\bi}$. In the second case $((X_{\et},\rho_{\et}),l)$ is mapped to $(X_{\et}\times X_{\m},(\rho_{\et},\rho_{\m}))$ with $X_{\m}=\hat{X}_{\et}$ and $\rho_{\m}=(p^l\cdot\rho_{\et})^{\vee}$. This morphism is a closed immersion. It is thus enough to show that each $k$-valued point of the left hand side is contained in the image. From the Hodge-Newton decomposition (see \cite{Katz}, Thm. 1.6.1) we obtain for each $ k $-valued point $(X,\rho)$ a decomposition $X=X_{\et}\times X_{\bi}\times X_{\m}$ and $\rho=\rho_{\et}\times \rho_{\bi}\times \rho_{\m}$ into the \'{e}tale, bi-infinitesimal, and multiplicative parts. The compatibility with the polarization then yields that up to some scalar $p^l$, the quasi- isogenies $\rho_{\m}$ and $\rho_{\et}$ are dual. From this the first isomorphism follows. The second isomorphism is shown by an easy calculation (compare \cite{modpdiv}, Section 3).
\end{proof}
The lemma reduces the questions after the global structure of $\mathcal{M}_{\red}$ to the same questions for $\mathcal{M}_{\mathbb{X}_{\bi},\red}$. Thus from now on we assume that $\mathbb{X}$ is bi-infinitesimal.


\section{The dense subscheme $\mathcal{S}_1$}\label{secs1}
In \cite{modpdiv}, 4.2 we define an open dense subscheme $\mathcal{S}_{\mathbb{X},1}^{\np}$ of $\mathcal{M}_{\mathbb{X},\red}^{\np}$. Let $\Lambda\subset N$ be the lattice associated to $x\in\mathcal{M}_{\mathbb{X},\red}^{\np}(k)$. Then $x\in\mathcal{S}_{\mathbb{X},1}^{\np}$ if and only if $a(\Lambda)=\dim_k(\Lambda/(F\Lambda+V\Lambda))=1$. As $F$ and $V$ are topologically nilpotent on $\Lambda$, this is equivalent to the existence of some $v\in \Lambda$ such that $\Lambda$ is the Dieudonn\'{e} submodule of $N$ generated by $v$. Note that $a(\Lambda)$ can also be defined as $\dim_k(\Hom(\alpha_p,X))$ where $X$ is the $p$-divisible group associated to $\Lambda$.

Let $$\mathcal{S}_1=\mathcal{S}_{\mathbb{X},1}^{\np}\cap\mathcal{M}_{\red}\subseteq\mathcal{M}_{\red}.$$ Then $\mathcal{S}_1$ is open in $\mathcal{M}_{\red}$.

\begin{lemma}\label{lemsdense}
The open subscheme $\mathcal{S}_1$ is dense in $\mathcal{M}_{\red}$. 
\end{lemma}
\begin{proof}
Recall that we assume $\mathbb{X}$ to be bi-infinitesimal. Let $(\overline{X},\overline{\rho})\in\mathcal{M}_{\red}(k)$ and let $\overline{\lambda}$ be a corresponding polarization of $\overline{X}$. Note that by \cite{deJong}, Lemma 2.4.4 there is an equivalence of categories between $p$-divisible groups over an adic, locally noetherian affine formal scheme $\Spf (A)$ and over $\Spec (A)$. From \cite{Oort1}, Corollary 3.11 we obtain a deformation $(X,\lambda)$ of $(\overline{X},\overline{\lambda})$ over $\Spec (k[[t]])$ such that the generic fiber satisfies $a=1$. It remains to show that we may also deform $\overline{\rho}$ to a quasi-isogeny $\rho$ between $(X,\lambda)$ and the constant $p$-divisible group $(\mathbb{X},\lambda_{\mathbb{X}})$ that is compatible with the polarizations. From \cite{OortZink}, Corollary 3.2 we obtain a deformation of $\overline{\rho}$ to a quasi-isogeny between $X$ and a constant $p$-divisible group after a base-change to the perfect hull of $k[[t]]$. After modifying this quasi-isogeny we may assume that the constant $p$-divisible group is the base change of $\mathbb{X}$. Let $x$ be the point of $\Spec(k[[t]]^{\rm perf})$ over the generic point of $\Spec(k[[t]])$. Then we may also modify the quasi-isogeny such that in $x$ it is compatible with the polarizations of the two groups in this point. Thus we obtain a $k[[t]]^{\rm perf}$-valued point of $\mathcal{M}_{\mathbb{X},\red}^{\np}$ such that the image of $x$ is in $\mathcal{M}_{\red}$. As $\mathcal{M}_{\red}$ is closed, this has to be a $k[[t]]^{\rm perf}$-valued point of $\mathcal{M}_{\red}$. If $\lambda_{\overline{X}},\lambda_{\overline{X}}'$ are two principal polarizations of $\overline{X}$, then $(\overline{X},\lambda_{\overline{X}})$ and $(\overline{X},\lambda_{\overline{X}}')$ differ by the isomorphism $(\lambda_{\overline{X}}')^{\vee}\circ\lambda_{\overline{X}}=\varphi$ with $\overline{\rho}^{-1}\varphi\overline{\rho}\in J$. Let $(\overline{X},\overline{\rho}')$ be the special point of the $k[[t]]^{\rm perf}$-valued point we constructed. Then there is an element of $J$ mapping $(\overline{X},\overline{\rho}')$ to $(\overline{X},\overline{\rho})$, where the compatibility with the polarizations follows from the previous argument. Thus we may assume that $\overline{\rho}'=\overline{\rho}$. Except at the special point, the $a$-invariant of the $p$-divisible group $X$ is 1. Thus this provides the desired deformation of $(\overline{X},\overline{\rho})$ to a point of $\mathcal{S}_1$. 
\end{proof}

To determine the dimension and the set of irreducible components of $\mathcal{M}_{\red}$ it is thus sufficient to consider $\mathcal{S}_1$. We proceed in the same way as for the moduli spaces of $p$-divisible groups without polarization. In contrast to the non-polarized case it turns out to be useful to use two slightly different systems of coordinates to prove the assertions on the dimension and on the set of irreducible components of $\mathcal{M}_{\red}$.

Let us briefly recall the main steps for the moduli spaces  $\mathcal{M}_{\mathbb{X},\red}^{\np}$ of non-polarized $p$-divisible groups. Their sets of irreducible components and their dimension are determined by studying $\mathcal{S}_{\mathbb{X},1}^{\np}$. In \cite{modpdiv}, 4 it is shown that the connected components of $\mathcal{S}_{\mathbb{X},1}^{\np}$ are irreducible and that $J^{\np}=\{j\in GL(N)|j\circ F=F\circ j\}$ acts transitively on them. The first step to prove this is to give a description of  $\mathcal{S}_{\mathbb{X},1}^{\np}(k)$. One uses that each such point corresponds to a lattice in $N$ with $a$-invariant $1$. As Dieudonn\'{e} modules, these lattices are generated by a single element and the description of the set of points is given by classifying these elements generating the lattices. The second step consists in the construction of a family in $\mathcal{M}^{\np}_{\mathbb{X},\red}$ to show that a set of points which seems to parametrize an irreducible component of $\mathcal{S}_{\mathbb{X},1}^{\np}$ indeed comes from an irreducible subscheme. More precisely, a slight reformulation of the results in \cite{modpdiv}, Section 4 yields the following proposition.

\begin{proposition}\label{propdisplay} 
Let $(N,F)$ be the isocrystal of a $p$-divisible group $\mathbb{X}$ over $k$. Let $m=v_p(\det F)$. Let $S=\Spec(R)\in \Nilp_W$ be a reduced affine scheme and $j\in J^{\np}$. Let $v\in N_R=N\otimes_{L}W(R)[\frac{1}{p}]$ such that in every $x\in S(k)$, the reduction $v_x$ of $v$ in $x$ satisfies that $$v_x\in j\Lambda_{\min}$$ and $$v_p(\det j)=\max\{v_p(\det j')\mid j'\in J^{\np} \text{ and }v_x\in j'\Lambda_{\min}\}.$$ Here, $\Lambda_{\min}\subset N$ is the lattice of the minimal $p$-divisible group in Remark \ref{remmin}. Let $\tilde{R}=\sigma^{-m}(R)$ be the unique reduced extension of $R$ such that $\sigma^m:\tilde{R}\rightarrow \tilde{R}$ has image $R$. Let $\tilde{v}\in N_{\tilde{R}}$ with $\sigma^m(\tilde{v})=v$. Then there is a morphism $\varphi:\Spec (\tilde{R})\rightarrow \mathcal{M}_{\mathbb{X},\red}^{\np}$ such that for every $x\in \Spec (\tilde{R})(k)$, the image $\varphi(x)$ correponds to the Dieudonn\'{e} module $\Lambda_x$ in $N$ generated by $v_x$. 
 
Assume in addition that $\mathbb{X}$ is principally polarized and that for every $x$, the Dieudonn\'{e} module $\Lambda_x$ corresponds to a point of $\mathcal{M}_{\red}$. Then $\varphi$ factors through $\mathcal{M}_{\red}$. 
\end{proposition}

Note that the second condition on $v_x$ (or more precisely the existence of the maximum) implies that the Dieudonn\'{e} submodule of $N$ generated by $v_x$ is a lattice.

\begin{proof}To prove the first assertion we may assume that $j=\id$. Note that $\tilde{v}$ satisfies the same conditions as $v$. The conditions on the $\tilde{v}_x$ are reformulated in \cite{modpdiv}, Lemma 4.7. The condition given there is exactly the condition needed in \cite{modpdiv}, 4.4 to construct a display over $S$ leading to the claimed morphism $\varphi$. It maps $x$ to the Dieudonn\'{e} lattice generated by $\sigma^m(\tilde{v}_x)=v_x$. The second assertion is trivial as $S$ is reduced and $\mathcal{M}_{\red}$ a closed subscheme of $\mathcal{M}_{\mathbb{X},\red}^{\np}$.
\end{proof}

\begin{remark}\label{remdualbasis}
We use the same notation as in the proposition. From \cite{modpdiv}, 4.4 we also obtain that under the conditions of Proposition \ref{propdisplay}, the elements $$v,Vv,\dotsc,V^{v_p(\det F)}v,Fv,\dotsc, F^{\dim N-v_p(\det F)-1}v$$ are a basis of the free $W(\tilde{R})[\frac{1}{p}]$-module $N_{\tilde{R}}$. They are the images of the standard basis of $N$ under some element of $GL(N_{\tilde{R}})$.

Let $v$ be as in Proposition \ref{propdisplay} with $N=N_0$. Then there are elements $y_i\in(N_1)_{\tilde{R}}$ which form a basis of $(N_1)_{\tilde{R}}$ which is dual to the basis $$(x_1,\dotsc,x_{\dim N_0})=(v,Vv,\dotsc,V^{v_p(\det F)}v,Fv,\dotsc, F^{\dim N_0-v_p(\det F)-1}v)$$ with respect to $\langle\cdot,\cdot \rangle$. In other words, the $y_i\in (N_1)_{\tilde{R}}$ are such that $\langle x_i,y_j\rangle=\delta_{ij}$.
\end{remark}


\section{Geometric points of $\mathcal{S}_1$}\label{secgeop}

\subsection{$N$ with an even number of supersingular summands} In this subsection we consider the case that $N$ has an even number of supersingular summands. 

Let $\Lambda\subset N$ be the lattice corresponding to a $k$-valued point of $\mathcal{M}_{\red}$. Then $\Lambda^{\vee}=c\Lambda$ for some $c\in L^{\times}$. Let $\Lambda_0=p_0(\Lambda)$ and $\Lambda_1=\Lambda\cap N_1$. For a subset $M$ of $N$ and $\delta\in\{0,1\}$ let 
\begin{equation}\label{gldual}
(M)_{\delta}^{\vee}=\{x\in N_{\delta}\mid \langle x,x'\rangle\in W\text{ for all }x'\in M\}.
\end{equation} 
Then $c\Lambda_1=(\Lambda_0)_1^{\vee}$. Hence $\Lambda_0$ and $\Lambda_1$ correspond to dual $p$-divisible groups, which implies $a(\Lambda_0)=a(\Lambda_1)$. 

The geometric points of $\mathcal{S}_1$ correspond to lattices $\Lambda$ which in addition satisfy $a(\Lambda)=1$. Especially, $a(\Lambda_0)=a(\Lambda_1)=1$. In this subsection we classify a slightly larger class of lattices. We fix a lattice $\Lambda_0\subset N_0$ with $a(\Lambda_0)=1$ and $c\in L^{\times}$. Then we consider all lattices $\Lambda\subset N$ with 
\begin{equation}\label{glreq}
p_0(\Lambda)=\Lambda_0 \text{ and }\Lambda^{\vee}=c\Lambda.
\end{equation}
Note that we have a description of the set of lattices $\Lambda_0\subset N_0$ with $a(\Lambda_0)=1$ from \cite{modpdiv}, see also Section \ref{secnp}.

The considerations above show that $\Lambda\cap N_1=\Lambda_1=c^{-1}(\Lambda_0)_1^{\vee}$ is determined by $\Lambda_0$ and $c$. Let $v_0$ be an element generating $\Lambda_0$ as a Dieudonn\'{e} module. If $v\in\Lambda$ with $p_0(v)=v_0$, then $\Lambda$ is generated by $v$ and $\Lambda_1$. Let $A$ be a generator of $\Ann(v_0)$ as in Lemma \ref{lemdefa}. We write $v=v_0+v_1$ for some $v_1\in N_1$. Then $Av=Av_1\in \Lambda_1$.
\begin{remark}
Let $\Lambda_0$, $c$, and $\Lambda_1$ be as above. Let $\Lambda$ be a Dieudonn\'{e} lattice with $p_0(\Lambda)=\Lambda_0$, $\Lambda\cap N_1\supseteq \Lambda_1$ and 
\begin{equation}\label{glbed2}
\Lambda^{\vee}\supseteq c\Lambda.\end{equation} The conditions imply that $\vol(\Lambda^{\vee})\leq\vol(c\Lambda)\leq \vol(c(\Lambda_0\oplus\Lambda_1))=\vol((\Lambda_0\oplus\Lambda_1)^{\vee})$. Dualizing the inequality for the first and last term, we see that all terms must be equal. Thus $\Lambda$ satisfies (\ref{glreq}) and $\Lambda\cap N_1=\Lambda_1$.
\end{remark}

The next step in the description of lattices with (\ref{glreq}) is to reformulate (\ref{glbed2}). It is equivalent to $\langle x,y\rangle\in c^{-1}W$ for all $x,y\in\Lambda$. Generators for $\Lambda$ as a $W$-module are given by $\Lambda_1$, the $F^iv$ with $0\leq i\leq n$, and the $V^iv$ with $0< i< m$ where $m$ is as in (\ref{gldefm}) and $n=h-m\geq m$. Thus (\ref{glbed2}) is equivalent to the condition that all scalar products between the various basis elements are in $c^{-1}W$. From the definition of $\Lambda_1$ we see that all products with elements of $\Lambda_1$ already satisfy this condition. By (\ref{glskpr}) it is enough to consider the products of $v$ with all other basis elements. Thus (\ref{glbed2}) is equivalent to
$$\langle v,F^iv\rangle\in c^{-1}W$$
and \begin{equation}\label{glf}\langle v,V^iv\rangle\in c^{-1}W\end{equation} for $n\geq i>0$. Furthermore, the equations for $V^i$ together with (\ref{glskpr}) imply those for $F^i$. 

All scalar products within one $N_i$ vanish. Hence the decomposition of $v$ together with (\ref{glskpr}) shows that (\ref{glf}) is equivalent to
\begin{equation}\label{glbed2b}
\langle v_0,V^iv_1\rangle-\langle V^iv_0,v_1\rangle\in c^{-1}W.
\end{equation}

An element $v_1\in N_1$ is uniquely determined by its scalar products with $F^iv_0$ for $i\in\{0,\dotsc,n-1\}$ and with $V^iv_0$ for $i\in \{1,\dotsc,m\}$. For $\phi\in\mathcal{D}$ let 
\begin{equation}\label{gldefxi}
\xi_{v_1}(\phi)=\langle \phi v_0, v_1\rangle.
\end{equation} 
Then $\xi_{v_1}$ is left-$W$-linear in $\phi$. The definition of $A$ implies that 
\begin{equation} \label{glxi2}
\xi_{v_1}(\psi A)=0
\end{equation}
for all $\psi\in\mathcal{D}$. We are looking for the set of $v_1$ satisfying (\ref{glbed2b}). In terms of $\xi_{v_1}$, this is
\begin{equation}
\label{glxi1}\xi_{v_1}(F^i)^{\sigma^{-i}}-\xi_{v_1}(V^i)\in c^{-1}W.
\end{equation}

\begin{lemma}\label{lemxiv}
\begin{enumerate}
\item Let $M$ be the set of $W$-linear functions $\xi:\mathcal{D}\rightarrow L$ with (\ref{glxi2}) and (\ref{glxi1}) for $i\leq n$. Then (\ref{gldefxi}) defines a bijection between $M$ and the set of elements $v_1\in N_1$ as above. 
\item Let $\overline{M}$ be the set of functions $\xi:\mathcal{D}\rightarrow L/c^{-1}W$ with the same properties as in (1). Then (\ref{gldefxi}) defines a bijection between $\overline{M}$ and the set of equivalence classes of elements $v_1$ as above. Here two such elements are called equivalent if their difference is in $c^{-1}(\Lambda_0)_1^{\vee}$.
\end{enumerate}
\end{lemma}
\begin{proof}
Let $\xi:\mathcal{D}\rightarrow L$ be given. An element $v_1$ of $N_1$ is uniquely determined by its scalar products with $v_0, Fv_0, \dotsc, F^{h-m-1}v_0$, $Vv_0,\dotsc, V^mv_0$. These $h$ values may be chosen arbitrarily. For the scalar products with the other elements of $\mathcal{D}v_0$, a complete set of relations is given by $\langle \psi Av_0,v_1\rangle=0$ for all $\psi\in\mathcal{D}$. This is equivalent to (\ref{glxi2}).
Furthermore, (\ref{glxi1}) is equivalent to the condition that the lattice generated by $\Lambda_1$ and $v_0+v_1$ satisfies  all required duality properties. 

To prove (2), we want to lift $\xi:\mathcal{D}\rightarrow L/c^{-1}W$ to a function with values in $L$. We lift the values of $\xi$ at $\phi\in\{V^m,V^{m-1},\dotsc,1,\dotsc,F^{h-m-1}\}$ arbitrarily. Then the lifts of the remaining values are uniquely determined by (\ref{glxi2}). As (\ref{glxi1}) was satisfied before, it still holds (as a relation modulo $c^{-1}W$) for the lifted functions. Then (1) implies the existence of $v_1$. Let now $w_1$ be a second element inducing $\xi$ (mod $c^{-1}W$). Then $\langle \phi v_0,w_1-v_1\rangle\in c^{-1}W$ for all $\phi\in \mathcal{D}$. Hence $w_1-v_1\in c^{-1}(\Lambda_0)_1^{\vee}$.
\end{proof}

\subsection{$N$ with an odd number of supersingular summands}\label{secodd}
As parts of this case are similar to the previous one, we mainly describe the differences.

We want to classify the lattices $\Lambda\subset N$ corresponding to $k$-valued points of $\mathcal{S}_1$. As before let $\Lambda_0=p_0(\Lambda)$ and $\Lambda_1=\Lambda\cap N_1$. Then $c\Lambda_1=(\Lambda_0)^{\vee}_1$. Besides,
\begin{equation}\label{glheihalb}
c\Lambda\cap N_{\frac{1}{2}}=(p_{\frac{1}{2}}(\Lambda))^{\vee}_{\frac{1}{2}}.
\end{equation} 
Here we use $(\cdot)^{\vee}_{\frac{1}{2}}$ analogously to (\ref{gldual}).

Again we already have a description of the Dieudonn\'{e} lattices $\Lambda_0\subset N_0$ with $a(\Lambda_0)=1$. We have to classify the $\Lambda$ correponding to some fixed $\Lambda_0$ and $c$, and begin by describing and normalizing the possible images under the projection to $N_0\oplus N_{\frac{1}{2}}$. Let $v\in\Lambda$ with $\mathcal{D}v=\Lambda$ and write $v=v_0+v_{\frac{1}{2}}+v_1$ with $v_i\in N_i$. Let $A$ with $\tilde{v}(A)= m$ be a generator of $\Ann(v_0)$ as in Lemma \ref{lemdefa}.
\begin{proposition}
\begin{enumerate}
\item There is a $j\in J$ such that $j(v)$ is of the form $v_0+\tilde{v}_{\frac{1}{2}}+\tilde{v}_1$ with $\tilde{v}_i\in N_i$ and $A\mathcal{D}\tilde{v}_{\frac{1}{2}}=\mathcal{D}A\tilde{v}_{\frac{1}{2}}$.
\item Let $j$ be as in the previous statement. Then $p_{\frac{1}{2}}(j\Lambda)$ is the unique Dieudonn\'{e} lattice in $N_{\frac{1}{2}}$ with $p_{\frac{1}{2}}(j\Lambda)^{\vee}=(cp^m)p_{\frac{1}{2}}(j\Lambda)$. Besides, $(j\Lambda)\cap N_{\frac{1}{2}}=p^mp_{\frac{1}{2}}(j\Lambda)$.
\end{enumerate}
\end{proposition}
\begin{proof}
Let $\tilde{v}_{\frac{1}{2}}\in N_{\frac{1}{2}}$ be such that $v'_{\frac{1}{2}}=v_{\frac{1}{2}}-\tilde{v}_{\frac{1}{2}}$ is in the kernel of $A$ and $A\mathcal{D}\tilde{v}_{\frac{1}{2}}=\mathcal{D}A\tilde{v}_{\frac{1}{2}}$. 

We first reduce the problem to the case where $N_0$ and $N_1$ are simple of slope $\frac{1}{2}$. Let $A_{\frac{1}{2}}$ be a generator of $\Ann(v'_{\frac{1}{2}})$ as in Lemma \ref{lemdefa}. As $A\in\Ann(v'_{\frac{1}{2}})$, we can write $A=\tilde{A}A_{\frac{1}{2}}$ with $\tilde{A}\in\mathcal{D}$. Then $\tilde{A}$ generates $\Ann(A_{\frac{1}{2}}v_0)$. We may write $v_0=v_0'+\tilde{v}_0$ with $\tilde{v}_0\in \tilde{N}_0$ and $A_{\frac{1}{2}}v_0'=0$. Then $v_0'$ generates a simple subisocrystal $N_0'$ of $N_0$ of slope $\frac{1}{2}$ and $N_0=N_0'\oplus\tilde{N}_0$. Let $N_1'$ be the subisocrystal of $N_1$ which is dual to $N_0'$. Then we want to show that the assertion of the proposition holds for some $j\in J\cap \End(N_0'\oplus N_{\frac{1}{2}}'\oplus N_1')$. To simplify the notation, we may assume that $N$ only consists of these three summands.

Recall that we can choose the following basis of $N$. Let $e_i^j\in N_j$ with $i=1,2$ and $j\in\{0,\frac{1}{2},1\}$. We may assume $F(e_i^j)=e_{i+1}^j$ using the convention $e_{i+2}^j=pe_i^j$. The basis can further be chosen such that the polarization is given by the property that all products of the generators vanish except $\langle e_i^j,e_{3-i}^{1-j}\rangle=1$ for all $j$. Replacing $j$ by $1-j$ induces an endomorphism of $N$ interchanging $N_0$ and $N_1$. Let $v_1'$ be the image of $v_0'$ under this isomorphism. Then $\Ann(v_1')=\Ann(v_{\frac{1}{2}}')=\Ann(v_0')$. Let $j^{-1}$ be given by mapping $v_0'$ to $v_0'+v_{\frac{1}{2}}'+av_1'$ where $a\in \mathbb{Q}_{p^2}$ is such that $\langle j^{-1}(v_0'),Fj^{-1}(v_0')\rangle=0$. Note that we can choose $a$ in this subfield of $L$ because $\Ann(v_0')=\Ann(v_{\frac{1}{2}}')=\Ann(v_1')$. Let further $j^{-1}\mid _{N_1}=\id$ and $j^{-1}(e_1^{\frac{1}{2}})=e_1^{\frac{1}{2}}+b$ where $b\in N_1$ is such that this element is orthogonal to the (simple) subisocrystal generated by $j^{-1}(v_0')$. As $a\in\mathbb{Q}_{p^2}$, we also have $b\in \mathbb{Q}_{p^2}e_1^1\oplus \mathbb{Q}_{p^2}e_2^1$. Thus $\Ann(j^{-1}e_1^{\frac{1}{2}})=\Ann(e_1^{\frac{1}{2}})=\mathcal{D}(F-V)$. For the remaining elements of $N$, we define $j^{-1}$ by requiring it to be linear and compatible with $F$. One then easily checks that we obtain a well-defined element of $J$. But $p_{\frac{1}{2}}(j(v_0+v_{\frac{1}{2}}+v_1))=v_{\frac{1}{2}}-v_{\frac{1}{2}}'=\tilde{v}_{\frac{1}{2}}$. Thus $j$ satisfies all properties of the first assertion.

For the second part of the proposition note that there is exactly one Dieudonn\'{e} lattice of each volume in $N_{\frac{1}{2}}$. Equivalently, for each $c$ there is exactly one $\Lambda\subset N_{\frac{1}{2}}$ with $\Lambda^{\vee}=c\Lambda$. Let $\Lambda_{\frac{1}{2}}\subseteq N_{\frac{1}{2}}$ be the lattice with $\Lambda_{\frac{1}{2}}^{\vee}=c\Lambda_{\frac{1}{2}}$. We have $\Lambda\cap N_{\frac{1}{2}}=c^{-1}(p_{\frac{1}{2}}(\Lambda))^{\vee}\subseteq p_{\frac{1}{2}}(\Lambda)$. This can be extended to a chain of inclusions $\Lambda\cap N_{\frac{1}{2}}\subseteq\Lambda_{\frac{1}{2}}\subseteq p_{\frac{1}{2}}(\Lambda)$ and the lengths of the two inclusions are equal. On the other hand, $p_{\frac{1}{2}}(\Lambda)$ contains $p_{\frac{1}{2}}(A\Lambda)=A(p_{\frac{1}{2}}(\Lambda))$. As $\tilde{v}(A)=m$, the length of this inclusion is $m$. Furthermore, $\Lambda\cap N_{\frac{1}{2}}=\Ann(Av_1)A\Lambda\subseteq \Ann(Av_1)Ap_{\frac{1}{2}}(\Lambda)$ is contained in $Ap_{\frac{1}{2}}(\Lambda)$. The length is at least $\tilde{v}(A_1)=m$ where $\Ann(Av_1)=\mathcal{D}A_1$. But on the other hand the duality relation for $\Lambda$ implies that $(p_{\frac{1}{2}}(A\Lambda))^\vee\supseteq cp_{\frac{1}{2}}(A\Lambda)$ or $p_{\frac{1}{2}}(A\Lambda)\subseteq \Lambda_{\frac{1}{2}}$. Comparing the two chains of inclusions, we see that $p_{\frac{1}{2}}(A\Lambda)=\Lambda_{\frac{1}{2}}$ and that the length of $\Lambda\cap N_{\frac{1}{2}}\subseteq Ap_{\frac{1}{2}}(\Lambda)$ is also $m$. This immediately implies the second part of the proposition.  
\end{proof}
For both Theorem \ref{thmirrkomp} and Theorem \ref{thmdim} it is enough to describe $\mathcal{S}_1$ up to the action of $J$. Thus we may assume that $j=1$ and that $v$ itself already satisfies the property of the proposition. Especially, $p_{\frac{1}{2}}(\Lambda)$ is then determined by $c$.

The element $v_{\frac{1}{2}}$ may be modified by arbitrary elements in $p_{\frac{1}{2}}(\Lambda\cap(N_{\frac{1}{2}}\oplus N_1))$ without changing $\Lambda$. Indeed, for each such element there is an element in $\Lambda$ whose projection to $N_0\oplus N_{\frac{1}{2}}$ is the given element. Thus for fixed $v_0$, the projection of $\Lambda$ to $N_0\oplus N_{\frac{1}{2}}$ is described by the element $v_{\frac{1}{2}}$ varying in the $W$-module $$p_{\frac{1}{2}}(\Lambda)/p_{\frac{1}{2}}(\Lambda\cap(N_{\frac{1}{2}}\oplus N_1))=p_{\frac{1}{2}}(\Lambda)/A(p_{\frac{1}{2}}(\Lambda))$$ of length $m$ which is independent of $\Lambda$. To choose coordinates for $v_{\frac{1}{2}}$ we use that it is isomorphic to $W/p^{\lfloor m/2\rfloor}W\oplus W/p^{\lceil m/2\rceil}W$. Under this isomorphism, the element $v_{\frac{1}{2}}$ is mapped to an element of the form $$\sum_{i=1}^{\lfloor m/2\rfloor}[y_i]p^{i-1}\oplus \sum_{\lfloor m/2\rfloor +1}^{m}[y_i]p^{i-\lfloor m/2\rfloor-1}.$$ Here we use that $k$ is perfect, and $[y_i]$ is the Teichm\"{u}ller representative of an element of $k$.

Note that $a(\Lambda)=1$ (or the condition that $j=1$) implies that $A(v_{\frac{1}{2}})$ is a generator of $\Lambda\cap N_{\frac{1}{2}}$ and not only an arbitrary element. This is an open condition on $p_{\frac{1}{2}}(\Lambda)/p_{\frac{1}{2}}(\Lambda\cap(N_{\frac{1}{2}}\oplus N_1)).$ More precisely, it excludes a finite number of hyperplanes (compare \cite{modpdiv}, Lemma 4.7).

Let now $v_{\frac{1}{2}}$ also be fixed. It remains to determine the set of possible $v_1$ such that $\Lambda=\mathcal{D}(v_0+v_{\frac{1}{2}}+v_1)$ is a lattice with $\Lambda^{\vee}=c\Lambda$. The same arguments as in the previous case show that $v_1$ can be chosen in an open subset of the set of $v_1$ with
\begin{equation}
\label{gl4.14}\langle v_0,\phi v_1\rangle +\langle v_1,\phi v_0\rangle\equiv- \langle v_{\frac{1}{2}},\phi v_{\frac{1}{2}}\rangle \pmod{c^{-1}W}.
\end{equation}
Let $\phi\in\mathcal{D}$ with $\tilde{v}(\phi)=2m$. Then $\phi v_{\frac{1}{2}}\in p^{ m}p_{\frac{1}{2}}(\Lambda)\subset c^{-1}\Lambda^{\vee}$. Especially, $\langle v_{\frac{1}{2}},\phi v_{\frac{1}{2}}\rangle$ is in $c^{-1}W$. This is later used in the form that $a_i=- \langle v_{\frac{1}{2}},F^i v_{\frac{1}{2}}\rangle$ satisfies (\ref{glcond1}).

Analogously to the previous case we use (\ref{gldefxi}) to define $\xi_v$. Then we also obtain the analogue of Lemma \ref{lemxiv}.


\section{The set of irreducible components}\label{secirr}

\begin{lemma}\label{lemaj}
Let $\Lambda\subset N_0\oplus N_1$ be a lattice generated by a sublattice $\Lambda_1\subset N_1$ and an element $v$ with $v=v_0+v_1$ for some $v_0\in N_0$ and $v_1\in N_1$. Let $\tilde{\Lambda}$ be generated by $\Lambda_1$ and $v_0+\tilde{v}_1$ for some $\tilde{v}_1\in N_1$. If $\xi_{\tilde{v}_1}(F^i)^{\sigma^{-i}}-\xi_{\tilde{v}_1}(V^i)=\xi_{v_1}(F^i)^{\sigma^{-i}}-\xi_{v_1}(V^i)$ for every $i\in \{1,\dotsc, h\}$ then there is a $j\in J$ with $j(\Lambda)=\tilde{\Lambda}$.
\end{lemma}
\begin{proof}
The assumption implies that $\langle v_0+\tilde{v}_1-v_1,\varphi(v_0+\tilde{v}_1-v_1)\rangle=0$ for $\varphi\in \{1,V,\dotsc, V^h\}$. By (\ref{glskpr}), the same holds for $\varphi\in \{F,\dotsc, F^h\}$. As $\dim N=2h$, the $\varphi(v_0+\tilde{v}_1-v_1)$ for these elements $\varphi\in\mathcal{D}$ generate $N'=\left(\mathcal{D}(v_0+\tilde{v}_1-v_1)\right)[1/p]\subseteq N$ as an $L$-vector space. Especially,  
\begin{equation}\label{glaniso}
\langle v_0+\tilde{v}_1-v_1,\varphi(v_0+\tilde{v}_1-v_1)\rangle=0
\end{equation}
for all $\varphi\in\mathcal{D}$. Let $A$ be a generator of $\Ann(v_0)$ as in Lemma \ref{lemdefa}. Then (\ref{glaniso}) for $\varphi=\varphi' A$ implies that $\langle v_0,\varphi' A(\tilde{v}_1-v_1)\rangle=0$ for all $\varphi'\in\mathcal{D}$. Thus $A(\tilde{v}_1-v_1)=0$.

Let $j\in GL(N)$ be defined by $v_0\mapsto v_0+\tilde{v}_1-v_1$, $j|_{N_1}=\id$, and $j\circ F=F\circ j$. To check that this is well-defined we have to verify that $Aj(v_0)=j(Av_0)=0$. But $A(j(v_0))=A(v_0+\tilde{v}_1-v_1)=0$. By definition $j$ commutes with $F$. Furthermore, (\ref{glaniso}) implies that $j\in G(L)$. Hence $j\in J$.
\end{proof}

For $v_1$ as above and $i\in\{1,\dotsc,h\}$ let 
\begin{equation}\label{gldefpsi}
\psi_i(v_1)=\xi_{v_1}(V^i)-\xi_{v_1}(F^i)^{\sigma^{-i}}.
\end{equation} 
Then the lemma can be reformulated as follows. Let $\Lambda$ and $\tilde{\Lambda}$ be two extensions of $\Lambda_0$ and $\Lambda_1$ as described in the previous section (or, in the case of an odd number of supersingular summands, two extensions of $\Lambda_0$ and $\Lambda_1$ associated to the same $v_{\frac{1}{2}}$) and let $v=v_0+v_1$ and $\tilde{v}=v_0+\tilde{v}_1$ (resp. $v=v_0+v_{\frac{1}{2}}+v_1$ and $\tilde{v}=v_0+v_{\frac{1}{2}}+\tilde{v}_1$) be the generators. Then $\psi_i(v_1)=\psi_i(\tilde{v}_1)$ for all $i$ implies that $\Lambda$ and $\tilde{\Lambda}$ are in one $J$-orbit. 

The next proposition implies that for each $(c_1,\dotsc,c_h)\in L^h$ there is a $v_1\in N_1$ with $\psi_i(v_1)=c_i$ for all $i$.

We fix an irreducible component of $\mathcal{S}_{\mathbb{X}_0,1}^{\np}$. Then \cite{modpdiv}, 4 describes a morphism from a complement of hyperplanes in an affine space to this irreducible component that is a bijection on $k$-valued points. Let $\Spec (R_0)$ be this open subscheme of the affine space. One first defines a suitable element $v_{0,R_0}\in N_0\otimes_{W}W(R_0)$. The morphism is then constructed in such a way that each $k$-valued point $x$ of $\Spec(R_0)$ is mapped to the lattice in $N_0$ generated by the reduction of $\sigma^m(v_{0,R_0})$ at $x$.
\begin{proposition}\label{lemirrkomp}
$R$ be a reduced $k$-algebra containing $\sigma^m(R_0)$. Let $c_1,\dotsc, c_h\in W(R)[1/p]$. Then there is an \'{e}tale extension $R'$ of $R$ and a $v_1\in N_{1,R'}$ with $\psi_i(v_1)=c_i$ for all $i$. Here, the $\psi_i$ are defined with respect to the universal element $\sigma^m(v_{0,R_0})\in (N_0)_{\sigma^m(R_0)}$. 
\end{proposition} 
For the proof we need the following lemma to simplify the occurring system of equations.
\begin{lemma}\label{lemglsys}
Let $m, n\in\mathbb{N}$ with $m\leq n$. For $0\leq i\leq m$ and $0\leq j\leq n$ let $P_{ij}(x)\in L[x]$ be a linear combination of the $\sigma^l(x)=x^{p^l}$ with $l\geq 0$. Assume that the coefficient of $x$ is zero for $j<i$ and in $W^{\times}$ for $i=j$. Consider the system of equations $$\sum_j P_{ij}(x_j)= a_i$$ with $a_i\in L$ and $i=0,\dotsc,i_0$ for some $i_0$. It is equivalent to a system of equations of the form $\sum_j Q_{ij}(x_j)= b_i$ with $b_i\in L$ such that the $Q_{ij}$ satisfy the same conditions as the $P_{ij}$ and in addition $Q_{ij}=0$ if $j<i$. 
\end{lemma}
\begin{proof}
We use a modification of the Gauss algorithm to show by induction on $\lambda$ that the system is equivalent to a system of relations of the form $\sum_j Q_{ij}^{\lambda}(x_j)= b_i^{\lambda}$ with $b_i^{\lambda}\in L$ such that the $Q_{ij}^{\lambda}$ satisfy the same conditions as the $P_{ij}$ and in addition $Q_{ij}^{\lambda}=0$ if $j<i$ and $j\leq \lambda$. For the induction step we have to carry out the following set of modifications for $j=\lambda+1$ and each $i> \lambda+1$. If $Q_{ij}^{\lambda}$ vanishes, we do not make any modification. We now assume $Q_{ij}^{\lambda}$ to be nontrivial. Let $\sigma^{l_i}(x)$ and $\sigma^{l_j}(x)$ be the highest powers of $x$ occurring in $Q_{ii}^{\lambda}$ and $Q_{ij}^{\lambda}$. If $l_i<l_j$, we modify the $j$th equation by a suitable multiple of $\sigma^{l_j-l_i}$ applied to the $i$th equation to lower $l_j$. Else we modify the $i$th equation by a suitable multiple of $\sigma^{l_i-l_j}$ applied to the $j$th equation to lower $l_i$. We proceed in this way as long as none of the two polynomials $Q_{ii}^{\lambda}$ and $Q_{ij}^{\lambda}$ becomes trivial. Note that the defining properties of the $P_{ij}$ are preserved by these modifications. As (by induction) $Q_{ij}^{\lambda}$ does not have a linear term, the linear term of $Q_{ii}^{\lambda}$ remains unchanged. Thus this process of modifications ends after a finite number of steps with equations $\sum_j Q_{ij}^{\lambda+1}(x_j)= b_i^{\lambda+1}$ which satisfy $Q_{ij}^{\lambda+1}=0$ for $j<i$ and $j\leq \lambda+1$. For $\lambda+1=n$, this is what we wanted.
\end{proof}

\begin{proof}[Proof of Proposition \ref{lemirrkomp}]
An element $v_1\in N_{1,R'}$ is determined by the values of $\xi_{v_1}$ at any $h$ consecutive elements of $\dotsc,F^2,F,1,V,V^2,\dotsc$. The other values of $\xi$ are then determined by $\xi_{v_1}(\phi A)=0$ for all $\phi\in\mathcal{D}$. Indeed, each of these equations for $\phi=F^i$ or $V^i$ for some $i$ gives a linear equation with coefficients in $L$ between the values of $\xi_{v_1}$ at $h+1$ consecutive elements of $\dotsc,F^2,F,1,V,V^2,\dotsc$. For this proof we want to take the values $\xi_{v_1}(F^i)$ for $i\in\{1,\dotsc,h\}$ as values determining $v_1$. Then all other values are linear combinations of these $\xi_{v_1}(F^i)$.
 
The definition of $\psi_{v_1}$ in (\ref{gldefpsi}) yields
\begin{equation*}
\xi_{v_1}(V^i)^{\sigma^i}=\xi_{v_1}(F^i)+\psi_i(v_1)^{\sigma^i}
\end{equation*} 
for $i\in\{1,\dotsc,h\}$. On the other hand, $\xi_{v_1}(V^i)^{\sigma^i}$ is a linear combination of the $\xi_{v_1}(F^j)^{\sigma^i}$ for $j\in \{1,\dotsc,h\}$. From this we obtain a system of $h$ equations for the $\xi_{v_1}(F^i)$ with $1\leq i\leq h$ of the same form as in Lemma \ref{lemglsys}. The resulting equations $\sum_j Q_{ij}(\xi_{v_1}(F^j))=b_i$ may be reformulated as $Q_{ii}(\xi_{v_1}(F^i))=c_i$ where $c_i$ also contains the summands corresponding to powers of $F$ larger than $i$. We can then consider these equations by decreasing induction on $i$. For each $i$, the polynomial $Q_{ii}(x)$ is a linear combination of powers of $x$ of the form $x^{\sigma^l}$, and its linear term does not vanish. Thus there is an \'{e}tale extension $R'$ of $R$ and $\xi_{v_1}(F^i)\in W(R')\otimes \mathbb{Q}$ with $v_p(\xi_{v_1}(F^i))\geq v_p(c_i)$ satisfying these equations. Note that $R'$ is in general an infinite extension of $R$, because the equations are between elements of $W(R)\otimes\mathbb{Q}$ and not over $R$ itself. Given $\xi_{v_1}$, Remark \ref{remdualbasis} shows that there is an element $v_1\in (N_1)_{R'}$ which induces $\xi_{v_1}$. Indeed, choose $v_1$ to be a suitable linear combination of the dual basis defined there.
\end{proof}

\subsection{Proof of Theorem \ref{thmirrkomp}}

To construct an irreducible subscheme of $\mathcal{S}_1$, we assume $c=1$. There is a $d\in\mathbb{N}$ such that for each $(c_1,\dotsc, c_h)\in (p^dW)^h$, the $v_1$ constructed in Proposition \ref{lemirrkomp} lies in the lattice $\Lambda_1\subset N_1$. Let $R_0$ as above. In the case of an even number of supersingular summands let $R_1=\sigma^m(R_0)$. Else let $(\sigma^m(v_{0,R_0}),v_{\frac{1}{2}})\in N_{ \sigma^m(R_0)[y_1,\dotsc,y_{ m}]}$ be the universal element (where $v_{\frac{1}{2}}$ is defined as in Section \ref{secodd}). The open condition on $\Spec (\sigma^m(R_0)[y_1,\dotsc,y_{ m}])$ that $Av_{\frac{1}{2}}$ is a generator and not only an element of $p_{\frac{1}{2}}(\Lambda\cap N_{\frac{1}{2}}\oplus N_1)$ is equivalent to the condition that $y_1$ does not lie in some finite-dimensional $\mathbb{F}_{p^2}$ sub-vector space of $k$ determined by the kernel of $A$. In this case let $R_1$ be the extension of $\sigma^m(R_0)$ corresponding to this affine open subscheme. Let in both cases $$R=R_1[x_{i,j}\mid i\in\{1,\dotsc,h\},j\in \{0,d-1\}].$$ Let $v_{0,R_0}$ as in Proposition \ref{lemirrkomp} and $c=1$. For $i\in\{1,\dotsc,h\}$ let $a_i=\sum_{j=0}^{d-1}[x_{i,j}]p^j$. Let $\Spec (R')$ and $v_1\in N_{R'}$ be as in the proposition. Let $v=\sigma^m(v_{0,R_0})+v_{1}$, resp. $v=\sigma^m(v_{0,R_0})+v_{\frac{1}{2}}+v_{1}$. Let $S=\Spec (R)$ be an irreducible component of the affine open subscheme of $\Spec (R')$ consisting of the points $x$ with $v_{1,x}\in(\mathcal{D}v_x)_1^{\vee}\setminus (F(\mathcal{D}v_x)_1^{\vee}+V(\mathcal{D}v_x)_1^{\vee})$. We denote the image of $v$ in $N_R$ also by $v$. As we already know that $\mathcal{S}_1$ is dense, this open subset is nonempty. Let $\tilde{R}$ be the inverse image of $R$ under $\sigma^h$ as in Proposition \ref{propdisplay}. Note that $v_p(\det F)=h$, whereas $v_p(\det F|_{N_0})=m$.

The next step is to define an associated morphism $\varphi:\tilde{R}\rightarrow \mathcal{M}_{\red}$ such that in each $k$-valued point $x$ of $S$, the image in $\mathcal{M}_{\red}(k)$ corresponds to the lattice generated by the reduction $v_x$ of $v$ at $x$. By Proposition \ref{propdisplay} it is enough to show that there is a $j\in J$ such that for each $x\in S(k)$, we have $v_x\in j \Lambda_{\min}$ and $v_p(\det j)=\max\{v_p(\det j')\mid v_x\in j'\Lambda_{\min}\}$. Let $\eta$ be the generic point of $S$ and let $j_{\eta}\in J$ be such a maximizing element for $\eta$. Then the same holds for each $k$-valued point in an open and thus dense subscheme of $S$. As the property $v_x\in j_{\eta}\Lambda_{\min}$ is closed, it is true for each $x\in S(k)$. In \cite{modpdiv}, 4 it is shown that for lattices $\Lambda\subset N$ with $a(\Lambda)=1$, the difference $\vol(\Lambda)-\max\{v_p(\det j')\mid \Lambda\subseteq j'\Lambda_{\min}\}$ is a constant only depending on $N$. In our case, the duality condition shows that $\vol (\mathcal{D} v_x)$ is constant on $S$ and only depending on $c$ and $N$. Thus the maximum is also constant. Hence in every $k$-valued point, $v_p(\det j_{\eta})$ is equal to this maximum, which is what we wanted for the existence of $\varphi:\Spec (\tilde{R})\rightarrow \mathcal{M}_{\red}$. We obtain an irreducible subscheme $\varphi(\Spec (\tilde{R}))$ of $\mathcal{S}_1\subseteq\mathcal{M}_{\red}$.

To show that $J$ acts transitively on the set of irreducible components we have to show that for each $x\in\mathcal{S}_1(k)$ there is an element $j\in J$ such that $jx$ lies in the image of $\varphi$. Let $\Lambda\subset N$ be the lattice corresponding to $x$. The first step is to show that there is a $j\in J$ such that $j(\Lambda)$ is selfdual (and not only up to a scalar $c(\Lambda)$). It is enough to show that there is a $j\in J$ such that $v_p(c(\Lambda))=v_p(c(j\Lambda))+1$. Such an element is given by taking the identity on $N_1$, multiplication by $p$ on $N_0$, and the map $e_i^j\mapsto e_{i+1}^j$ on $N_{\frac{1}{2}}$. Here we use the notation of Remark \ref{remmin} for the basis of $N$. Next we want to modify $\Lambda_0$. We have $a(p_0(\Lambda))=a(\Lambda\cap N_1)=1$. From the classification of lattices with $a=1$ in \cite{modpdiv}, 4 we obtain that $J_{\mathbb{X}_0}^{\np}$ (which may be considered as a subgroup of $J$) is acting transitively on the set of irreducible components of $\mathcal{M}^{\np}_{\mathbb{X}_0,\red}$. Thus by multiplying with such an element we may assume that $\Lambda_0$ lies in the fixed irreducible component chosen for Proposition \ref{lemirrkomp}. Recall from Section \ref{secodd} that in the case of an odd number of supersingular summands, there is a $j\in J$ mapping the element $v_{\frac{1}{2}}$ to the irreducible family described there. It remains to study the possible extensions between the lattices $\Lambda_0$ and $\Lambda_1$ (or in the second case between the sublattice of $N_0\oplus N_{\frac{1}{2}}$ determined by $\Lambda_0$ and $v_{\frac{1}{2}}$ and $\Lambda_1$). They are given by the associated elements $v_1$. Fix a generating element $\sigma^m(v_0)$ of $\Lambda_0$ (in the second case also an element $v_{\frac{1}{2}}$) and let $v_1$ be an element associated to an extension $\Lambda$ with $a(\Lambda)=1$. Then it is an immediate consequence of Lemma \ref{lemaj} and the construction of $S$ that there is an element of $J$ mapping $\Lambda$ to a lattice associated to a point of $S$ inducing the same $\psi_i$ as $\Lambda$. Thus the image of $S$ under $J$ is $\mathcal{S}_1$, which proves the theorem.\qed


\section{Dimension}
To determine the dimension of $\mathcal{S}_1$ and of $\mathcal{M}_{\red}$ we have to classify the elements $v_1$ of Section \ref{secgeop} up to elements in $c^{-1}\Lambda_1$ and not up to the (locally finite) action of $J$ which we used in Section \ref{secirr}. To do so, it is more useful to use the values of $\xi_{v_1}$ as coordinates instead of the values of $\psi_{v_1}$.

We investigate the set of possible values $\xi(\phi)\in L/c^{-1}W$ for $\phi\in\mathcal{D}$ using decreasing induction on $\tilde{v}(\phi)\geq 0$. Here, $\tilde{v}$ is as in (\ref{gldefvtilde}). Recall from Lemma \ref{lemxiv} (2) that the use of functions $\xi$ with values in $L/c^{-1}W$ instead of $L$ corresponds to considering $v_1$ as an element of $N_1/\Lambda_1$. But as $v_1$ and $v_1+\delta$ with $\delta\in\Lambda_1$ lead to the same lattice $\Lambda$, this is sufficient to determine the set of possible extensions of $\Lambda_0$ and $\Lambda_1$. 

Instead of equations (\ref{glxi2}) and (\ref{glxi1}) we consider the following slightly more general problem to treat at the same time the case of an odd number of supersingular summands. There, (\ref{glxi1}) is replaced by (\ref{gl4.14}). We want to consider $W$-linear functions $\xi:\mathcal{D}\rightarrow L/c^{-1}W$ with
\begin{eqnarray}
\label{glxiodd1}\xi(F^i)-\xi(V^i)^{\sigma^i}&\equiv& a_i\pmod{ c^{-1}W}\\
\label{glxiodd2}\xi(\psi A)&\equiv& 0\pmod{ c^{-1}W}
\end{eqnarray}
for all $\psi\in\mathcal{D}$. Here $a_i\in L$ are such that 
\begin{eqnarray}
\label{glcond1}
a_ip^{j_i}\in c^{-1}W &\text{ if }&2j_i+i\geq 2 m.
\end{eqnarray}

Let $\mathcal{D}^i=\{\phi\in\mathcal{D}\mid \tilde{v}(\phi)\geq i\}$. We call a $W$-linear function 
$$\xi^{i_0}:\mathcal{D}^{i_0}\rightarrow L/(c^{-1}W)$$ satisfying (\ref{glxiodd1}) and (\ref{glxiodd2}) a partial solution of level $i_0$. Then the induction step consists in determining the possible partial solutions $\xi^{i_0}$ of level $i_0$ leading to a fixed solution of level $i_0+1$. Note that the assumption on $a_i$ implies that there exists the trivial partial solution $\xi^{2 m}\equiv 0$ of level $2 m$ inducing partial solutions of all higher levels. Recall that we assumed $F$ and $V$ to be elementwise topologically nilpotent on $N$. Thus for each function $\xi$ with (\ref{glxiodd1}) and (\ref{glxiodd2}) there is a level $i$ such that $\xi$ induces the trivial partial solution of level $i$.

Assume that we already know the $\xi(\phi)$ for $\tilde{v}(\phi)>i_0$ and want to determine its possible values for $\tilde{v}(\phi)=i_0$. Then we know in particular $\xi(p\phi)=p\xi(\phi)\in L/c^{-1}W$, or $\xi(\phi)\in L/p^{-1}c^{-1}W$. We want to determine the possible liftings modulo $c^{-1}W$.

A basis of the $ k $-vector space $\mathcal{D}^{i_0}/\mathcal{D}^{i_0+1}$ is given by the $i_0+1$ monomials $$F^{i_0},pF^{i_0-2},\dotsc, pV^{i_0-2},V^{i_0}.$$ Equation (\ref{glxiodd1}) leads to $\lfloor i_0/2\rfloor$ relations between the values of $\xi$ on these monomials. Recall that $\tilde{v}(A)= m$. Thus if $\tilde{v}(\phi)= i_0-m$ for some $\phi \in\mathcal{D}$, (\ref{glxiodd2}) leads to a relation between the value of $\xi$ on $\LT(\phi A)\in\mathcal{D}^{i_0}$ and values on $\mathcal{D}^{i_0+1}$. As the $\xi$ are linear, it is sufficient to consider the $\max\{0,i_0- m+1\}$ relations for $\phi\in\{F^{i_0-m},pF^{ i_0-m-2},\dotsc,V^{ i_0-m}\}\cap\mathcal{D}^{i_0-m}$. This count of relations leads to the notation $$r(i_0)=\lfloor i_0/2\rfloor+\max\{0,i_0- m+1\}.$$ Then $i_0+1\leq r(i_0)$ is equivalent to $i_0\geq 2 m$.

The following proposition is the main tool to prove Theorem \ref{thmdim} on the dimension of the moduli spaces.
\begin{proposition}\label{propmain}
\begin{enumerate}
\item Let $i_0\geq 2 m$. Then there is a partial solution $\xi^{i_0}$ of (\ref{glxiodd1}) and (\ref{glxiodd2}) of level $i_0$. If we fix $\xi^{i_0}$ and an $l\in\mathbb{N}$ with $l\geq i_0$, there are only finitely many other partial solutions $\tilde\xi^{i_0}$ of level $i_0$ such that $\xi^{i_0}-\tilde\xi^{i_0}$ induces the trivial partial solution of level $l$ of the associated homogenous system of equations.
\item Let $i_0+1>r(i_0)$ and let $\xi^{i_0+1}$ be a partial solution of (\ref{glxiodd1}) and (\ref{glxiodd2}) of level $i_0+1$. Then to obtain a partial solution $\xi^{i_0}$ of level $i_0$ inducing $\xi^{i_0+1}$, one may choose the lifts to $L/c^{-1}W$ of the values of $\xi^{i_0}$ at the first $i_0+1-r(i_0)$ monomials $p^{\alpha} V^{\beta}$ with $2\alpha+\beta=i_0$ and $\beta\leq 2(i_0-r(i_0))+1$ arbitrarily. Each of the remaining values lies in some finite nonempty set depending on the values on the previous monomials. 
\end{enumerate}
\end{proposition}
\begin{proof}
Note that the existence statement in the first assertion is satisfied as the condition on the $a_i$ yields that there is the trivial solution of level $2 m$. We show the two assertions simultaneously. Let $\xi^{i_0+1}$ be a fixed partial solution of level $i_0+1$ for any $i_0$. It is enough to show that for a lift $\xi^{i_0}$, the values of the first $\max \{0,i_0+1-r(i_0)\}$ variables can be chosen arbitrarily, and that the remaining values then lie in some finite set. If $i_0+1>r(i_0)$, we have to show that this finite set is nonempty. We investigate the relations (\ref{glxiodd1}) and (\ref{glxiodd2}) more closely. The first set of relations shows that $\xi^{i_0}(p^aF^b)$ with $2a+b=i_0$ is determined by $\xi^{i_0}(p^aV^b)$. Thus it is sufficient to consider this latter set of values.
Besides, we have to consider (\ref{glxiodd2}) for $\psi\in\{V^{i_0- m},pV^{i_0- m-2},\dotsc,F^{i_0- m}\}$. For $B\in\mathcal{D}$ let $\LT(B)$ as in Lemma \ref{lemdefa}. Then the equations for the values of $\xi^{i_0}$ relate $\xi^{i_0}(\LT(\psi A))$ to something which is known by the induction hypothesis. Recall the description of $\LT(A)$ from Lemma \ref{lemdefa} (3). Let $h'$ be the number of supersingular summands of $N_0$. Let $j\geq 0$ with $i_0- m-j\geq 0$. Then $\LT(V^{i_0- m-j}F^{j}A)$ is a linear combination of $V^{i_0-j}F^{j},\dotsc,V^{i_0-j-h'}F^{j+h'}$ whose coefficients are Teichm\"{u}ller representatives of elements of $k$. Furthermore, the coefficients of $\xi^{i_0}(V^{i_0-j}F^{j})$ and $\xi^{i_0}(V^{i_0-j-h'}F^{j+h'})$ are units in $W$. Using (\ref{glxiodd1}) we may replace values of $\xi^{i_0}$ at monomials in $F$ by $\sigma$-powers of the values of the corresponding monomials in $V$. We thus obtain a relation between a polynomial in the remaining $\lceil (i_0+1)/2\rceil$ values of $\xi^{i_0}$ and an expression which is known by induction. For $2j\leq i_0$, the first summand $\xi^{i_0}(V^{i_0-j}F^j)$ remains the variable associated to the highest power of $V$ which occurs linearly in this polynomial. In the following we ignore all equations for $2j>i_0$. They only occur for $i_0> 2 m$, a case where we only want to prove the finiteness of the set of solutions. The system of equations with $2j\leq i_0$ is of the form requested in Lemma \ref{lemglsys}. The proof of this Lemma for coefficients in $L/c^{-1}W$ is the same as for coefficients in $L$. Thus we obtain that the lifts of the values at the $i_0+1-r(i_0)$ variables associated to the largest values of $j$ can be chosen freely and the other ones have to satisfy some relation of the form $Q_{ii}(x)\equiv b_i$ for some given $b_i$. As the $Q_{ii}$ have a linear term they are nontrivial. This implies that the set of solutions of these equations is nonempty and finite. 
\end{proof}

\subsection{Proof of Theorem \ref{thmdim}}
By Lemma \ref{lemsdense} it is enough to show that $\mathcal{S}_1$ is equidimensional of the claimed dimension. From \cite{modpdiv}, 4 we obtain that the connected components of $\mathcal{S}_{\mathbb{X},1}^{\np}$ are irreducible. The discrete invariant with values in $J^{\np}/(J^{\np}\cap \Stab(\Lambda_{\min}))$ distinguishing the components is given by $$\Lambda\mapsto j_{\Lambda} \text{ with } \Lambda\subseteq j_{\Lambda}\Lambda_{\min} \text{ and } v_p(\det j_{\Lambda})=\max \{v_p(\det j)\mid j\in J^{\np}, \Lambda\subseteq j\Lambda_{\min}\}.$$ Especially, $j_{\Lambda}$ is constant on each connected component of $\mathcal{S}_1\subseteq \mathcal{S}_1^{\np}$. Besides, $p_0(j_{\Lambda}\Lambda_{\min})$ determines the connected component of $p_0(\Lambda)=\Lambda_0$ inside $\mathcal{S}_{\mathbb{X}_0,1}$. Thus we may fix an irreducible component of $\mathcal{S}_{\mathbb{X}_0,1}$ and determine the dimension of the union of connected components of $\mathcal{S}
_1$ such that $\Lambda_0$ is in this fixed component. Let $R_0$ and $R_1$, $v_0$ and $v_{\frac{1}{2}}$ be as in the proof of Theorem \ref{thmirrkomp}. Again we use the functions $\xi$ defined with respect to $\sigma^m(v_0)$ instead of $v_0$. Fix an arbitrary partial solution $\xi^{2m}$ of level $2m$. Let $$R_2=R_1[x_{i\beta}\mid i\geq 0,1\leq \beta \leq i+1-r(i)].$$ We use decreasing induction on $i$ to lift $\xi^{2m}$ to a partial solution of level $i$ over some extension $R_2^{i}$ of $R_2$. Let $R_2^{2m}=R_2$. Assume that a lift $\xi^{i+1}$ is given. Then Proposition \ref{propmain} shows that the values at $i+1-r(i)$ monomials with $\tilde{v}=i$ may be lifted arbitrarily to a value of $\xi^{i}$. If $p^{\alpha}V^{\beta}$ with $2\alpha+\beta=i$ and $\beta \leq i+1-r(i)]$ is such a monomial we write (using the induction hypothesis) $\xi^{i+1}(p^{\alpha+1}V^{\beta})=\sum_{i< v_p(c^{-1})}[a_i]p^i$ with $a_i\in R_{2}^{i+1}$. Then we choose 
$$\xi^i(p^{\alpha}V^{\beta})=\sum_{i< v_p(c^{-1})}[a_i]p^{i-1}+[x_{i\beta}]p^{v_p(c^{-1})-1}.$$
Let $R_2^{i}$ be the extension of $R_2^{i+1}$ given by adjoining further variables $x_{i\beta}$ for larger $\beta$ parametrizing the other values of the lift of $\xi^{i+1}$ to $\xi^{i}$ and with relations as in Proposition \ref{propmain}(2) and its proof. The special form of the $Q_{ii}$ implies that $R_2^i$ is a finite \'{e}tale extension of $R_2^{i+1}$. Let $R_3=R_2^{0}$. Let $v_{1,R_3}\in N_{1,R_3}$ be such that $\xi_{v_{1,R_3}}=\xi^0$. Its existence follows again from the existence of the dual basis in Remark \ref{remdualbasis}. Let $v=\sigma^m(v_{0,R_0})+v_{1,R_3}$, or $v=\sigma^m(v_{0,R_0})+v_{\frac{1}{2}}+v_{1,R_3}$. As in the proof of Theorem \ref{thmirrkomp} let $S=\Spec(R)$ be an irreducible component of the affine open subscheme of $\Spec (R_3)$ over which $\mathcal{D}v$ contains $(\mathcal{D}v)_1^{\vee}$. We only have to consider the case that this subscheme is nonempty because we already know that $\mathcal{S}_1$ itself is dense and thus nonempty. Let $\tilde{R}=\sigma^{-m}(R)$ as in Proposition \ref{propdisplay}. The same argument as in the proof of Theorem \ref{thmdim} shows that there is a morphism $\varphi:\Spec{(\tilde{R})}\rightarrow \mathcal{M}_{\red}$ mapping $x\in \Spec{(\tilde{R})}(k)$ to the lattice generated by $v_x$. The finiteness statements in Proposition \ref{propmain} imply that for each given $y\in\mathcal{S}_1$ (and thus given $\xi$) there is an open neighborhood in $\mathcal{S}_1$ which only contains points of $\varphi(\Spec{(\tilde{R})})$ for a finite number of choices of $S$. Besides, the construction of $R_3$ together with the description of the $k$-valued points of $\mathcal{S}_1$ shows that for each $y\in\mathcal{S}_1(k)$ there is exactly one choice of $S$ and one point $x\in \Spec{(\tilde{R})}(k)$ such that $\varphi(x)=y$. Thus $\dim \mathcal{M}_{\red}=\dim \mathcal{S}_1$ is the maximum of $\dim \Spec{(\tilde{R})}$ for all irreducible components $S$. It remains to show that this is equal to the right hand side of (\ref{gldimform}). Note that $R_i$ is equidimensional for $i=0,1,2,3$. From the construction of $S$ we see that in case of an even number of supersingular summands, 
\begin{equation}\label{glm-1}
\dim \Spec{(\tilde{R})}=\dim S=\dim R_3=\dim R_2=\dim \mathcal{M}_{\mathbb{X}_0,\red}^{\np}+\sum_{i\geq 0}\max\{0,i+1-r(i)\}.
\end{equation}
In the other case,
\begin{equation}\label{glm-2}
\dim \Spec{(\tilde{R})}=\dim S=\dim R_3=\dim R_2=\dim \mathcal{M}_{\mathbb{X}_0,\red}^{\np}+\sum_{i\geq 0}\max\{0,i+1-r(i)\}+m. 
\end{equation}
The last summand corresponds to the choice of $v_{\frac{1}{2}}$.
 
From the decomposition of $N$ into $l$ simple summands $N^j$ we obtain a decomposition $N_0=\bigoplus_{j=1}^{l_0}N^j$ with $l_0=\lfloor l/2\rfloor$. Let again $\lambda_j=m_j/(m_j+n_j)$ be the slope of $N^j$. Recall from \cite{modpdiv}, Theorem B that $$\dim \mathcal{M}_{\mathbb{X}_0,\red}^{\np}=\sum_{j=1}^{l_0}\frac{(m_j-1)(n_j-1)}{2}+\sum_{\{j,j'\mid j< j'\leq l_0\}}m_jn_{j'}.$$ We denote the right hand side of (\ref{gldimform}) by $D$. Let us first consider the case of an even number of supersingular summands. Then by the symmetry of the Newton polygon we obtain
\begin{align*} 
D-\dim \mathcal{M}_{\mathbb{X}_0,\red}^{\np}&=\frac{1}{2}\sum_{j< j'\leq l}m_jn_{j'}+\frac{ m}{2}-\sum_{j< j'\leq l_0}m_jn_{j'}.
\intertext{ Again by the symmetry of the Newton polygon this is equal to}
&=\sum_{j=1}^{l_0}\sum_{j'=l_0+1}^{l}\frac{m_jn_{j'}}{2}+\frac{ m}{2}\\
&=\sum_{j,j'=1}^{l_0}\frac{m_jm_{j'}}{2}+\frac{ m}{2}\\
&=\frac{m(m+1)}{2}.
\end{align*}
In the other case, the same calculation shows that
\begin{align*} 
D-\dim \mathcal{M}_{\mathbb{X}_0,\red}^{\np}
&=\frac{ m( m+1)}{2}+2\sum_{j=1}^{l_0}\frac{m_jn_{l_0+1}}{2}\\
&=\frac{ m( m+1)}{2}+ m.
\end{align*}
In the last step we used that $N^{l_0+1}$ is supersingular, hence $n_{l_0+1}=1$. 

On the other hand (\ref{glm-1}) implies that in the case of an even number of supersingular summands,
\begin{eqnarray*}
\dim \mathcal{M}_{\red}-\dim\mathcal{M}_{\mathbb{X}_0,\red}^{\np}&=&\sum_{i\geq 0}\max\{0,i+1-r(i)\}\\
&=&\sum_{i=0}^{m-1}\left(\left\lfloor \frac{i}{2}\right\rfloor+1\right)+ \sum_{i=m}^{2m-1}\left(\left\lfloor \frac{i}{2}\right\rfloor-i+m\right)\\
&=&m+\sum_{i=0}^{2m-1}\left\lfloor \frac{i}{2}\right\rfloor-\sum_{i=0}^{m-1}i\\
&=&\frac{m(m+1)}{2}.
\end{eqnarray*}
The same calculation with (\ref{glm-2}) shows that for an odd number of supersingular summands $$\dim \mathcal{M}_{\red}-\dim\mathcal{M}_{\mathbb{X}_0,\red}^{\np}=\frac{ m( m+1)}{2}+ m.$$ Together with the calculation of $D-\dim \mathcal{M}_{\mathbb{X}_0,\red}^{\np}$, this implies Theorem \ref{thmdim}.\qed


\section{Connected components}

In this section we determine the set of connected components of $\mathcal{M}_{\red}$. Again, the essential tools are the description of $\mathcal{S}_1$ and the comparison with the non-polarized case. The reduction to the bi-infinitesimal case in Section \ref{secredbi} shows that Theorem \ref{thmconncomp} follows from the next theorem.

\begin{theorem}\label{thmconncomp2}
Let $\mathbb{X}$ be bi-infinitesimal and non-trivial. Then 
\begin{eqnarray*}
\kappa:\mathcal{M}_{\red}( k )&\rightarrow& \mathbb{Z}\\
\Lambda&\mapsto&v_p(c(\Lambda)),
\end{eqnarray*} where $\Lambda^{\vee}=c(\Lambda)\cdot\Lambda$, induces a bijection $$\pi_0(\mathcal{M}_{\red})\cong \mathbb{Z}.$$
\end{theorem}

\begin{proof}[Proof of Theorem \ref{thmconncomp2}]
From Theorem \ref{thmirrkomp} we obtain a $J$-invariant surjection $\pi:J\twoheadrightarrow \pi_0(\mathcal{M}_{\red})$. Besides, the map $\kappa$ induces a surjection $\pi_0(\mathcal{M}_{\red})\rightarrow \mathbb{Z}$. We choose the basis of $N$ described in Remark \ref{remmin} to identify $N$ with $L^{2h}$. Then the standard lattice is the lattice $\Lambda_{\min}$ corresponding to a minimal $p$-divisible group and to a point of $\mathcal{M}_{\red}$. An element $j\Lambda_{\min}$ with $j\in J$ is in the kernel of $\kappa$ if and only if $(j\Lambda_{\min})^{\vee}=j\Lambda_{\min}$. This is equivalent to $j\Lambda_{\min}=j'\Lambda_{\min}$ for some $j'\in J\cap Sp_{2h}(L)$. Thus we have to show that $J\cap Sp_{2h}(L)$ is mapped to a single connected component of $\mathcal{M}_{\red}$. There is the following explicit description of $J\cap Sp_{2h}(L)$. Let $N=\bigoplus_{\lambda} N(\lambda)$ be the decomposition of $N$ into isoclinic subisocrystals $N(\lambda)$ of slope $\lambda$. For $\lambda<\frac{1}{2}$ let $J_{\lambda}=GL_{d_{\lambda}}(D_{\lambda})$. Let $J_{\frac{1}{2}}=Sp_{d_{\frac{1}{2}}}(D_{\frac{1}{2}})$. In both cases $D_{\lambda}$ is a central simple algebra over $\mathbb{Q}_p$ of invariant $\lambda$ and $d_{\lambda}$ is the number of summands in a decomposition of $N(\lambda)$ into simple isocrystals. Then $\prod_{\lambda\leq \frac{1}{2}}J_{\lambda}\rightarrow J\cap Sp_{2h}(L)$, which maps $((j_{\lambda})_{\lambda<\frac{1}{2}},j_{\frac{1}{2}})$ to $((j_{\lambda}),j_{\frac{1}{2}},(j_{\lambda})^{\vee})$, is an isomorphism.

Let $M\subset Sp_{2h}(L)\cong Sp(N)$ be the Levi subgroup $GL(N_0)$ or $GL(N_0)\times Sp(N_{\frac{1}{2}})$ with associated standard parabolic subgroup $P$. As a Borel subgroup we choose the one fixing the flag associated to the ordered isotropic subset $e_{1}^1,\dotsc,e^1_{h_1},e_1^2,\dotsc,e^{l_0}_{h_{l_0}}$ (and in the case of an odd number of supersingular summands the additional element $e_{1}^{l_0+1}$) of $N$. Here we use the notation $l_0=\lfloor \frac{l}{2}\rfloor$ and the basis of Remark \ref{remmin}. Let $U$ be the unipotent radical of $P$. We denote the simple subgroup $GL(N_0)\times \id$ of $M$ by $M_0$. Note that if the number of supersingular summands of $N$ is even, $M=M_0$ is the Siegel Levi subgroup. The explicit description of $J$ shows that $J\cap Sp_{2h}(L)$ is generated by $J\cap M$, $J\cap K$ and $J\cap U$. As $\pi$ is $J$-equivariant, it is enough to show that each of these three subgroups lies in a single connected component.

Our choice of the basis of $N$ implies that the stabilizer of the base point of $\mathcal{M}_{\red}$ is $\Stab (\Lambda_{\min})=K$. Thus the surjection $\pi$ maps $J\cap K$ to the component of the identity.

Next we consider $J\cap M$. Note that in the second case, $(J\cap M)/K=(J\cap M_0)/K$. Thus it is sufficient to consider $J\cap M_0$. We use the canonical inclusion $\mathcal{M}^{\np}_{\mathbb{X}_0,\red}\hookrightarrow \mathcal{M}_{\red}$ of (\ref{glincx0}) whose image contains $J\cap M_0$. By the results of \cite{modpdiv}, Section 3, $\pi_0(\mathcal{M}^{\np}_{\mathbb{X}_0,\red})$ is isomorphic to $\mathbb{Z}$ or trivial. The latter case occurs only if $M_0$ is trivial. The discrete invariant in the former case is given by the valuation of the determinant of the corresponding element of $GL(N_0)$ or equivalently by the volume of the corresponding lattice in $N_0$. To show that $J\cap M_0$ lies in a single connected component in $\mathcal{M}_{\red}$ also in the case of nontrivial $M_0$, we have to connect two elements representing $0$ and $1$ in $\mathbb{Z}$ by a family. To do so, we may assume that $N_0$ and $N_1$ are simple isocrystals of slopes $\lambda$ and $1-\lambda$ for some $\lambda$. Indeed, in the general case we can construct a connecting family by taking the family in the special case and adding constant trivial families in the other simple summands of $N$. Write $\lambda=\frac{m}{h'}$ with $(m,h')=1$. (If $M=M_0$, we have $h=h'$, otherwise $h'=h-1$.) We choose the same basis of $N_{\lambda}$ and $N_{1-\lambda}$ as before. Then $e^1_1,\dotsc e_{h'}^1$ form the basis of $N_0$ and $e^2_1,\dotsc, e^2_{h'}$ the basis of $N_1$ with $\langle e^1_i, e^2_j\rangle=\delta_{i,h'-j}$ for all $i,j\in\{1,\dotsc,h'\}.$ Let $\Lambda_0$ be the lattice generated by these basis elements, and let $\Lambda_1$ be generated by $e^1_2,\dotsc,e^1_{h'},e^1_{h'+1},e^2_0,e^2_1,\dotsc,e^2_{h'-1}$. One easily checks that both lattices are selfdual. We have to show that $\Lambda_0$ and $\Lambda_1$ are in the same connected component of $\mathcal{M}_{\red}$. In \cite{modpdiv}, Lemma 3.4 these two lattices are connected by an affine line in $\mathcal{M}_{\mathbb{X},\red}^{\np}$ mapping $x\in\mathbb{A}^1(k)$ to the lattice $\Lambda_x$ generated by $\Lambda_0\cap \Lambda_1$ and $(1-x)e^1_1+xe^2_0$. It is thus enough to show that the image of this map $\mathbb{A}^1\rightarrow \mathcal{M}_{\mathbb{X},\red}^{\np}$ lies in $\mathcal{M}_{\red}$. To check this we have to show that the $\Lambda_x$ are selfdual, which is an easy consequence of the definition of the two bases. Thus $J\cap M$ lies in the connected component of the identity in $\mathcal{M}_{\red}$.

For $a\in\mathbb{Z}$ let $j'_a\in J\cap Z(M)$ where $Z(M)$ is the center of $M$ with $j'_a=p^a$ in the case of an even number of supersingular summands of $N$ and $j'_a=(p^a,\id)$ in the case of an odd number of supersingular summands. Let $j\in J\cap U$. Then there is an $a$ such that $j'_aj(j'_a)^{-1}\in K$. Especially it is in the connected component of $\id\in\mathcal{M}_{\red}$. As $j'_a$ is in $K$, we obtain that $j$ is also in this connected component.   

Together we obtain that $J\cap M$, $J\cap K$ and $J\cap U$ are all in the kernel of $\pi$. Thus the $J$-invariant surjection $\pi$ induces a bijection $\pi_1(G)\cong \pi_0(\mathcal{M}_{\red})$. 
\end{proof}

\end{document}